\documentstyle[12pt]{amsart}
\def\NZQ{\Bbb}               
\def\NN{{\NZQ N}}

\def\ZZ{{\NZQ Z}}
\def\RR{{\NZQ R}}

\def\opn#1#2{\def#1{\operatorname{#2}}} 
\opn\chara{char}
\opn\length{\ell}
\opn\pd{pd}
\opn\rk{rk}
\opn\projdim{proj\,dim}
\opn\rank{rank}
\opn\depth{depth}
\opn\grade{grade}
\opn\height{height}
\opn\embdim{emb\,dim}
\opn\codim{codim}

\opn\Tr{Tr}
\opn\bigrank{big\,rank}
\opn\superheight{superheight}\opn\lcm{lcm}
\opn\trdeg{tr\,deg}%
\opn\reg{reg}
\opn\lreg{lreg}
\opn\div{div}
\opn\Div{Div}
\opn\cl{cl}
\opn\Cl{Cl}
\opn\Spec{Spec}
\opn\Supp{Supp}
\opn\supp{supp}
\opn\Sing{Sing}
\opn\Ass{Ass}
\opn\Ann{Ann}
\opn\Rad{Rad}
\opn\Soc{Soc}
\opn\Ker{Ker}
\opn\Coker{Coker}
\opn\Im{Im}
\opn\Hom{Hom}
\opn\Tor{Tor}
\opn\Ext{Ext}
\opn\End{End}
\opn\Aut{Aut}
\opn\id{id}

\opn\nat{nat}
\opn\pff{pf}
\opn\Pf{Pf}
\opn\GL{GL}
\opn\SL{SL}
\opn\mod{mod}
\opn\ord{ord}
\opn\aff{aff}
\opn\con{conv}
\opn\relint{relint}
\opn\st{st}
\opn\lk{lk}
\opn\cn{cn}
\opn\core{core}
\opn\vol{vol}
\opn\link{link}
\opn\star{star}
\opn\gr{gr}

\def\pot#1#2{#1[\kern-0.28ex[#2]\kern-0.28ex]}

\opn\dirlim{\underrightarrow{\lim}}
\opn\inivlim{\underleftarrow{\lim}}
\let\sect=\cap

\let\mcone= * 
\let\to=\rightarrow

\def\Implies{\ifmmode\Longrightarrow \else
     \unskip${}\Longrightarrow{}$\ignorespaces\fi}
\def\implies{\ifmmode\Rightarrow \else
     \unskip${}\Rightarrow{}$\ignorespaces\fi}
\def\iff{\ifmmode\Longleftrightarrow \else
     \unskip${}\Longleftrightarrow{}$\ignorespaces\fi}

\let\:=\colon
\newtheorem{Theorem}{Theorem}[section]
\newtheorem{Lemma}[Theorem]{Lemma}
\newtheorem{Corollary}[Theorem]{Corollary}
\newtheorem{Proposition}[Theorem]{Proposition}
\newtheorem{Remark}[Theorem]{Remark}

\newtheorem{Example}[Theorem]{Example}
\newtheorem{Examples}[Theorem]{Examples}
\newtheorem{Definition}[Theorem]{Definition}

\let\epsilon\varepsilon
\let\phi=\varphi
\let\kappa=\varkappa
\def\qed{\ifhmode\textqed\fi
   \ifmmode\ifinner\quad\qedsymbol\else\dispqed\fi\fi}
\def\textqed{\unskip\nobreak\penalty50
    \hskip2em\hbox{}\nobreak\hfil\qedsymbol
    \parfillskip=0pt \finalhyphendemerits=0}
\def\dispqed{\rlap{\qquad\qedsymbol}}

\opn\inii{in}
\opn\inim{inm}
\opn\set{set}
\opn\sort{sort}
\opn\conv{conv}
\opn\dis{dis}

\begin{document}

\title{Discrete Polymatroids}
\author{J\"urgen Herzog and Takayuki Hibi}
\address{J\"urgen Herzog, Fachbereich Mathematik und Informatik,  
Universit\"at-GHS Essen, 45117 Essen, Germany}
\email{juergen.herzog@@uni-essen.de}
\address{Takayuki Hibi, Department of Mathematics, Graduate School of Science,
Osaka University, Toyonaka, Osaka 560-0043, Japan}
\email{hibi@@math.sci.osaka-u.ac.jp}

\maketitle
\begin{abstract}
The discrete polymatroid is a multiset analogue
of the matroid.  Based on the polyhedral theory
on integral polymatroids
developed in late 1960's and in early 1970's,
in the present paper the
combinatorics and algebra on discrete polymatroids
will be studied.
\end{abstract}


\section*{Introduction}
Matroid theory
is one of the most fascinating research areas
in combinatorics.
Let $[n] = \{1, 2, \ldots, n \}$ and $2^{[n]}$
the set of all subsets of $[n]$.
For $A \subset [n]$ write $|A|$ for the cardinality
of $A$.  A {\em matroid} on the ground set $[n]$ is a nonempty subset
${\cal M} \subset 2^{[n]}$ satisfying
\begin{enumerate}
\item[(M1)]
if $F_1 \in {\cal M}$ and $F_2 \subset F_1$,
then $F_2 \in {\cal M}$;
\item[(M2)]
if $F_1$ and $F_2$ belong to ${\cal M}$ and
$|F_1| < |F_2|$, then there is
$x \in F_2 \setminus F_1$ such that
$F_1 \cup \{ x \} \in {\cal M}$.
\end{enumerate}
The members of ${\cal M}$ are the
{\em independent sets} of ${\cal M}$.
A {\em base} of ${\cal M}$ is a maximal
independent set of ${\cal M}$.
It follows from (M2) that if $B_1$ and $B_2$
are bases of ${\cal M}$, then
$|B_1| = |B_2|$.
The set of bases of ${\cal M}$
possesses the ``exchange property'' as follows:
\begin{enumerate}
\item[(B)]
If $B_1$ and $B_2$ are bases of ${\cal M}$
and if $x \in B_1 \setminus B_2$, then
there is $y \in B_2 \setminus B_1$
such that
$(B_1 \setminus \{ x \}) \cup \{ y \}$
is a base of ${\cal M}$.
\end{enumerate}
Moreover, the set of bases of ${\cal M}$
possesses the ``symmetric exchange property''
as follows:
\begin{enumerate}
\item[(S)]
If $B_1$ and $B_2$ are bases of ${\cal M}$
and if $x \in B_1 \setminus B_2$, then
there is $y \in B_2 \setminus B_1$
such that both
$(B_1 \setminus \{ x \}) \cup \{ y \}$
and
$(B_2 \setminus \{ y \}) \cup \{ x \}$
are bases of ${\cal M}$.
\end{enumerate}
On the other hand,
given a nonempty set
${\cal B} \subset 2^{[n]}$,
there exists a matroid ${\cal M}$ on $[n]$
with ${\cal B}$ its set of bases
if and only if ${\cal B}$ possesses
the exchange property (B).

Let $\varepsilon_1, \ldots, \varepsilon_n$
denote the canonical basis vectors of
$\RR^n$.  If we associate each $F \subset [n]$
with $\sum_{i \in F} \varepsilon_i$,
then a matroid on $[n]$ may be regarded as
a set of $(0, 1)$-vectors of $\RR^n$.
Considering nonnegative integer vectors
instead of $(0, 1)$-vectors enables us
to define the concept of
discrete polymatroids.
Let $\RR_+^n$ denote the set of vectors
$u = (u_1, \ldots, u_n) \in \RR^n$
with each $u_i \geq 0$. If $u =(u_1, \ldots, u_n)$
and $v = (v_1, \ldots, v_n)$
are two vectors belonging to $\RR_+^n$, then 
we write $u \leq v$ if all components $v_i - u_i$
of $v - u$ are nonnegative and, moreover,
write $u < v$ if $u \leq v$ and $u \neq v$.
We say that $u$ is a {\em subvector} of $v$
if $u \leq v$.
The {\em modulus} of
$u = (u_1, \ldots, u_n) \in \RR_+^n$
is $|u| = u_1 + \cdots + u_n$.
Also, let $\ZZ_+^n = \RR_+^n \cap \ZZ^n$. 

A {\em discrete polymatroid} on the ground set $[n]$
is a nonempty finite set $P \subset \ZZ_+^n$
satisfying
\begin{enumerate}
\item[(D1)]
if $u  \in P$ and
$v  \in \ZZ_+^n$ with $v \leq u$, then
$v \in P$;
\item[(D2)]
if $u = (u_1, \ldots, u_n) \in P$ and
$v = (v_1, \ldots, v_n) \in P$ with
$|u| < |v|$, then there is $i \in [n]$
with $u_i < v_i$ such that
$u + \varepsilon_i \in P$.
\end{enumerate}
A {\em base} of $P$ is a vector $u \in P$
such that $u < v$ for no $v \in P$.
It follows from (D2) that if
$u_1$ and $u_2$ are bases of $P$,
then $|u_1| = |u_2|$.
Discrete polymatroids
can be also characterized in terms of
the exchange property of their bases.
In fact,
a nonempty finite set $B \subset \ZZ_+^n$
is the set of bases of a discrete
polymatroid on $[n]$ if and only if
$B$ satisfies
\begin{enumerate}
\item[(i)]
all $u \in B$ have the same modulus;
\item[(ii)]
if $u = (u_1, \ldots, u_n) \in P$ and
$v = (v_1, \ldots, v_n) \in P$ belong to
$B$ with $u_i > v_i$, then there is
$j \in [n]$ with $u_j < v_j$
such that
$u - \varepsilon_i + \varepsilon_j \in B$.
\end{enumerate}
Consult Theorem 2.3 for the details.

In, e.g., \cite[Chapter 18]{Welsh}
and \cite[p. 382]{Oxley}
we can find the concept of polymatroids,
which originated in Edmonds \cite{Edmonds}.
A {\em polymatroid} on $[n]$ is
a nonempty  compact subset ${\cal P} \subset \RR_+^n$
satisfying
\begin{enumerate}
\item[(P1)]
if $u \in {\cal P}$
and
$v \in \RR_+^n$
with  $v \leq u$, then
$v \in {\cal P}$;
\item[(P2)]
if $u = (u_1, \ldots, u_n) \in {\cal P}$ and
$v = (v_1, \ldots, v_n) \in {\cal P}$ with
$|u| < |v|$, then there is $i \in [n]$
with $u_i < v_i$ and
$0 < N < v_i - u_i$
such that
$u + N \varepsilon_i \in {\cal P}$.
\end{enumerate}
A brief summary of fundamental materials
on polymatroids is given in Section $1$.
It follows that a polymatroid ${\cal P}$
on $[n]$ is a convex polytope in $\RR^n$.
We say that a polymatroid
${\cal P} \subset \RR_+^n$
is {\em integral}
if all vertices of ${\cal P}$
are integer  vectors of $\RR^n$.

Now, a combinatorial aspect of the present
paper is to understand the relation
between discrete polymatroids and
integral polymatroids.
Our first fundamental theorem is

\medskip

\noindent
{\bf Theorem 3.4.}
{\em A nonempty finite set
$P \subset \ZZ_{+}^{n}$ is a discrete polymatroid
if and only if $\conv(P) \subset \RR_{+}^{n}$
is an integral polymatroid with
$\conv(P) \cap \ZZ^{n} = P$.
Here $\conv(P)$ is the convex hull of
$P$ in $\RR^n$.}

\medskip

Moreover, Theorem 3.4 together with
combinatorics on integral polymatorids
yields our second fundamental theorem,
the symmetric exchange theorem for
discrete polymatroids, as follows:

\medskip

\noindent
{\bf Theorem 4.1.}
{\em
If $u = (u_1, \ldots, u_n)$ and
$v = (v_1, \ldots, v_n)$ are bases of a discrete
polymatroid $P \subset \ZZ_+^n$, then
for each $i \in [n]$ with $u_i > v_i$
there is $j \in [n]$ with $u_j < v_j$
such that both $u - \varepsilon_i + \varepsilon_j$
and $v - \varepsilon_j + \varepsilon_i$
are bases of $P$.
}

\medskip

On the other hand, an algebraic aspect of
the present paper is to study Ehrhart rings
and base rings together with toric ideals
of discrete polymatroids.
Let $K$ be a field and let $t_1, \ldots, t_n$
and $s$ indeterminates over $K$.
If $u = (u_1, \ldots, u_n) \in \ZZ_+^n$,
then $t^u = t_1^{u_1} \cdots t_n^{u_n}$.
Let $P \subset \ZZ_+^n$ be a discrete
polymatroid and $B$ its set of bases.
Let $K[P]$ denote the homogeneous
semigroup ring $K[t^u s \: u \in P]$.
The {\em base ring} of $P$ is the subalgebra
$K[B] = K[t^u s \: u \in B]$
$( \cong K[t^u \: u \in B])$
of $K[P]$.
Let $K[x_u \: u \in B]$
denote the polynomial ring
in $|B|$ variables over $K$
and write
$I_B \subset K[x_u \: u \in B]$
for the {\em toric ideal}
of $K[B]$, i.e., $I_B$ is the kernel of
the surjective $K$-algebra
homomorphism $\xi \:
K[x_u \: u \in B]
\to
K[B]$ defined by
$\xi(x_u) = t^u$
for all $u \in B$.
Clearly,
all symmetric relations belong to $I_B$.
More precisely,
if $u = (u_1, \ldots, u_n)$ and
$v = (v_1, \ldots, v_n)$ belong to $B$
with $u_i > v_i$ and $u_j < v_j$
and if
both $u - \varepsilon_i + \varepsilon_j$
and $v - \varepsilon_j + \varepsilon_i$
belong to $B$, then
the quadratic binomial
$x_u x_v -
x_{u - \varepsilon_i + \varepsilon_j}
x_{v - \varepsilon_j + \varepsilon_i}$
belong to $I_B$.
White \cite{White} conjectured that
for a matroid all symmetric exchange relations
generates $I_B$.  It is natural to conjecture
that this is also holds for discrete polymatroids.
Much stronger,
commutative algebraists cannot escape
from the temptation to study the following problems:
\begin{enumerate}
\item[(a)]
Does the toric ideal $I_B$ of a discrete polymatroid
possess a quadratic Gr\"obner basis?
\item[(b)]
Is the base ring $K[B]$ of a discrete
polymatroid Koszul?
\end{enumerate}
These must be the most attractive research
problems on base rings of discrete polymatroids.
One of our results on toric ideals
of discrete polymatroids is

\medskip

\noindent
{\bf Theorem 5.3.}
{\em
{\em (a)}
Suppose that each matroid has the property that
the toric ideal of its base ring is generated
by symmetric exchange relations, then this is
also true for each discrete polymatroid.

{\em (b)}
If $P \subset \ZZ_+^n$ is a discrete polymatroid
whose set of
bases $B$ satisfies the strong exchange property
(i.e.,
if $u = (u_1, \ldots, u_n),
v = (v_1, \ldots, v_n) \in B$,
then for all $i$ and $j$
with $u_i > v_i$ and $u_j < v_j$
one has 
$u - \varepsilon_i + \varepsilon_j\in B$), then
\begin{enumerate}
\item[(i)]
$I_B$ has a quadratic Gr\"obner basis,
and hence $K[B]$ is Koszul;
\item[(ii)]
$I_B$ is generated by symmetric exchange relations.
\end{enumerate}
}
\medskip

It is know that if
${\cal P}_1, \ldots, {\cal P}_k
\subset \RR_+^n$
are polymatroid then their polymatroid sum
\[
{\cal P}_1 \vee \cdots \vee {\cal P}_k
=
\{ x = \sum_{i=1}^k x_i \:
x_i \in {\cal P}_i, 1 \leq i \leq k \}
\]
is again a polymatroid.  Moreover,
if each ${\cal P}_i$ is integral, then
${\cal P}_1 \vee \cdots \vee {\cal P}_k$
is also integral and
for each integer vector
$x \in {\cal P}_1 \vee \cdots \vee {\cal P}_k$
there exist integer vectors $x_i \in {\cal P}_i$,
$1 \leq i \leq k$, with $x = \sum_{i=1}^k x_i$.
An immediate consequence of this fact is
the following important

\medskip

\noindent
{\bf Theorem 6.1.}
{\em
If $P \subset \ZZ_+^N$ is a discrete polymatroid,
then the homogeneous semigroup ring $K[P]$ is normal.
}

\medskip

\noindent
{\bf Corollary 6.2.}
{\em
The base ring $K[B]$ of a discrete polymatroid $P$
is normal.
}

\medskip

Thus in particular both $K[P]$ and $K[B]$
are Cohen--Macaulay.
We will also compare algebraic properties
of $K[P]$ with those of $K[B]$.
See Theorem 6.3 for the details.

Once we know that both $K[P]$ and $K[B]$
are Cohen--Macaulay,
it is natural to study the problem
when $K[P]$ is Gorenstein and when
$K[B]$ is Gorenstein.
A combinatorial criterion for $K[P]$ to be
Gorenstein  is obtained in Theorem 7.3.
To classify all Gorenstein base rings seems,
however, quite difficult.  We will
introduce the concept of ``generic''
discrete polymatroids and find
a combinatorial characterization
of discrete polymatroids $P$
satisfying the condition that
(i) $P$ is generic and (ii)
the base ring of $P$ is Gorenstein.
See Theorem 7.6.

Finally, in Section $8$,
two simple techniques to construct
discrete polymatroids will be studied.

\section{Polymatroids}
The present section is a brief summary,
based on \cite[Chapter 18]{Welsh},
of fundamental materials on polymatroids.
Fix an integer $n > 0$ and set
$[n] = \{ 1, 2, \ldots, n\}$.
The canonical basis vectors of $\RR^n$
will be denoted by
$\varepsilon_1, \ldots, \varepsilon_n$.
Let $\RR_+^n$ denote the set of those vectors
$u =(u_1, \ldots, u_n) \in \RR^n$ with each
$u_i \geq 0$,
and $\ZZ_+^n = \RR_+^n \cap \ZZ^n$.
For a vector
$u =(u_1, \ldots, u_n) \in \RR_+^n$
and for a subset $A \subset [n]$,
we set
\[
u(A) = \sum_{i \in A} u_i.
\]
Thus in particular $u(\{i\})$,
or simply $u(i)$, is the $i$th component
$u_i$ of $u$.
The {\em modulus} of $u$ is
\[
|u| = u([n]) = \sum_{i=1}^n u_i.
\]
Let $u =(u_1, \ldots, u_n)$
and $v = (v_1, \ldots, v_n)$
be two vectors belonging to $\RR_+^n$.
We write $u \leq v$ if all components $v_i - u_i$
of $v - u$ are nonnegative and, moreover,
write $u < v$ if $u \leq v$ and $u \neq v$.
We say that $u$ is a {\em subvector} of $v$
if $u \leq v$.
In addition, we set
\begin{eqnarray*}
u \vee v & = &
(\max\{u_1,v_1\},\ldots, \max\{u_n,v_n\}), \\
u \wedge v & = &
(\min\{u_1,v_1\},\ldots, \min\{u_n,v_n\}).
\end{eqnarray*}
Thus $u \wedge v \leq u \leq u \vee v$
and $u \wedge v \leq v \leq u \vee v$.

\begin{Definition}
\label{defpolymat}
{\em A {\em polymatroid} on the {\em ground set} $[n]$
is a nonempty compact subset ${\cal P}$ in $\RR_+^n$,
the set of {\em independent vectors}, such that
\begin{enumerate}
\item[(P1)]
every subvector of an independent vector is
independent;
\item[(P2)]
if $u, v \in {\cal P}$ with $|v| > |u|$, then
there is a vector $w \in P$ such that
\[
u < w \leq u \vee v.
\]
\end{enumerate}}
\end{Definition}

A {\em base} of a polymatroid ${\cal P} \subset \RR_+^n$
is a maximal independent vector of ${\cal P}$,
i.e., an independent vector $u \in {\cal P}$ with
$u < v$ for no $v \in {\cal P}$.
Every base of ${\cal P}$ has the same modulus
$\rank({\cal P})$, the {\em rank} of ${\cal P}$.
In fact, if $u$ and $v$ are bases of ${\cal P}$
with $|u| < |v|$,
then by (P2) there exists $w \in P$ with
$u < w \leq u\vee v$,
contradicting the maximality of $u$.

Let ${\cal P} \subset \RR_+^n$ be a polymatroid
on the ground set $[n]$.
Let $2^{[n]}$ denote the set of all subsets of $[n]$.
The {\em ground set rank function} of ${\cal P}$
is a function $\rho : 2^{[n]} \to \RR_{+}$
defined by setting
\[
\rho(A) = \max\{ v(A) \: v \in {\cal P} \}
\]
for all $\emptyset \neq A \subset [n]$ together with
$\rho(\emptyset) = 0$.

\begin{Proposition}
\label{submodular}
{\em (a)}
Let ${\cal P} \subset \RR_+^n$ be a polymatroid
on the ground set $[n]$
and $\rho$ its ground set rank function.  Then $\rho$
is nondecreasing, i.e., if $A \subset B \subset [n]$,
then $\rho(A) \leq \rho(B)$, and is submodular, i.e.,
\[
\rho(A) + \rho(B) \geq \rho(A \cup B) + \rho(A \cap B)
\]
for all $A, B \subset [n]$.
Moreover, ${\cal P}$ coincides with the compact set
\begin{eqnarray}
\{ x \in \RR_+^n \:
x(A) \leq \rho(A), A \subset [n] \}.
\end{eqnarray}

{\em (b)}
Conversely, given a nondecreasing and submodular
function
$\rho : 2^{[n]} \to \RR_{+}$
with $\rho(\emptyset) = 0$,
the compact set {\em (1)} is a polymatroid
on the ground set
$[n]$ with $\rho$ its ground set rank function.
\end{Proposition}

It follows from Proposition \ref{submodular} (a) that
a polymatroid ${\cal P} \subset \RR_+^n$
on the ground set $[n]$ is a convex
polytope in $\RR^n$.
In addition, the set of bases of ${\cal P}$
is a face of ${\cal P}$ with  supporting hyperplane
\[
\{ x =(x_1, \ldots, x_n) \in \RR^n \:
\sum_{i=1}^n x_i = \rank({\cal P})\}.
\]
We refer the reader to, e.g., \cite{Hibi} for basic
terminologies on convex polytopes.

How can we find the vertices of a polymatroid?
A complete answer was obtained by Edmonds \cite{Edmonds}.
We will associate any permutation
$\pi = (i_1, \ldots, i_n)$ of $[n]$
with
$A_\pi^1 = \{ i_1 \},
A_\pi^2 = \{ i_1, i_2 \}, \ldots,
A_\pi^n = \{ i_1, \ldots, i_n \}$.

\begin{Proposition}
\label{vertex}
Let ${\cal P} \subset \RR_+^n$ be a polymatroid
on the ground set
$[n]$ and $\rho$ its ground set rank function.
Then the vertices of ${\cal P}$ are all points
$v = v(k,\pi) \in \RR_+^n$, where
$v = (v_1, \ldots, v_n)$ and
\begin{eqnarray*}
v_{i_1} & = & \rho(A_\pi^1), \\
v_{i_2} & = & \rho(A_\pi^2) \, - \,
\rho(A_\pi^1), \\
v_{i_3} & = & \rho(A_\pi^3) \, - \,
\rho(A_\pi^2), \\
& & \cdots \\
v_{i_k} & = & \rho(A_\pi^k) \, - \,
\rho(A_\pi^{k-1}), \\
v_{i_{k+1}} = v_{i_{k+2}} = & \cdots &
= v_{i_n} = 0,
\end{eqnarray*}
and $k$ ranges over the integers belonging to $[n]$,
and
$\pi = (i_1, \ldots, i_n)$ ranges over all
permutations of $[n]$.
In particular the vertices of the face
of ${\cal P}$
consisting of all bases of ${\cal P}$ are
all points $v = v(n,\pi) \in \RR_+^n$, where
$\pi$ ranges over all
permutations of $[n]$.
\end{Proposition}

We say that a polymatroid is {\em integral}
if all of its vertices have
integer coordinates; in other words,
a polymatroid is integral if and only if
its ground set rank function is integer valued.

Let ${\cal P}_1, \ldots, {\cal P}_k$ be polymatroids
on the ground set $[n]$.  The {\em polymatroid sum}
${\cal P}_1 \vee \cdots \vee {\cal P}_k$
of ${\cal P}_1, \ldots, {\cal P}_k$
is the compact subset in $\RR_+^n$ consisting of
all vectors $x \in \RR_+^n$ of the form
\[
x = \sum_{i=1}^k x_i, \, \, \, \, \, \, \, \, \, \,
x_i \in {\cal P}_i.
\]

\begin{Proposition}
\label{sum}
Let ${\cal P}_1, \ldots, {\cal P}_k$ be polymatroids
on the ground set $[n]$ and
$\rho_i$ the ground set rank function
of ${\cal P}_i$, $1 \leq i \leq k$.  Then
the polymatroid sum
${\cal P}_1 \vee \cdots \vee {\cal P}_k$
is a polymatroid on $[n]$ and
$\sum_{i=1}^k \rho_i$
is its ground set rank function.
Moreover, if each ${\cal P}_i$ is integral, then
${\cal P}_1 \vee \cdots \vee {\cal P}_k$ is integral,
and for each integer vector
$x \in {\cal P}_1 \vee \cdots \vee {\cal P}_k$
there exist integer vectors $x_i \in {\cal P}_i$,
$1 \leq i \leq k$, with $x = \sum_{i=1}^k x_i$.
\end{Proposition}
\section{Discrete  polymatroids}

In this section we introduce discrete polymaroids. They may be viewed as generalizations of 
matroids. 

\begin{Definition}
{\em Let $P$ be a nonempty finite set of integer  vectors in $\RR^n_+$ which contains with 
each $u\in P$ all its integral subvectors. The set $P$ is called a {\em discrete polymatroid} 
on the ground set $[n]$ if for all $u,v\in P$ with $|v|>|u|$, there is a vector $w\in P$ such 
that 
\[
u<w\leq u\vee v.
\]}
\end{Definition}

A {\em base} of $P$ is a vector $u \in P$
such that $u < v$ for no $v \in P$.
We  denote  $B(P)$ the set of bases of a discrete polymatroid $P$. Just as for polymatroids 
one proves that all  elements in $B(P)$ have the same modulus. This common number is called 
the {\em rank of $P$}.  

The following simple lemma is useful for induction arguments.

\begin{Lemma}
\label{cut}
Let $P$ be a  discrete  polymatroid.
\begin{enumerate}
\item[(a)]  Let $d\leq \rank P$. Then the set $P'=\{u\in P\: |u|\leq d\}$ is a discrete 
polymatroid of rank $d$ with the set of bases $\{u\in P\: |u|=d\}$.
\item[(b)] Let $x\in P$. Then the set $P_x=\{v-x \: v\geq x\}$ is a discrete polymatroid 
of rank $d-|x|$.
\end{enumerate}
\end{Lemma}

\begin{pf}
(a) Let $u,v\in P$ with $d\geq |v|>|u|$. There exists $w\in P$ such that $u<w\leq u\wedge v$. 
Since $w>u$, and since $P$ contains all subvectors of $w$, there exists an integer $i$
such that $u+\epsilon_i\leq w$.  Then $u<u+\epsilon_i\leq u\wedge v$, and since 
$u+\epsilon_i\leq d$, it belongs to $P'$. This proves that  $P'$ is a discrete polymatroid. It 
is clear that $\{u\in P\: |u|=d\}$ is the set of bases of $P'$. 

(b) Let $u',v'\in P_x$ with $|v'|>|u'|$, and set $u=u'+x$ and $v=v'+x$. Then $u,v\in P$ and 
$|v|>|u|$. Hence there exists $w\in P$ with $u<w\leq u\wedge v$. Set $w'=w-x$. Then $w'\in 
P_x$ and $u'<w'\leq u'\wedge v'$. 
\end{pf}

Discrete polymatroids can be characterized in terms of their set of bases as follows:

\begin{Theorem}
\label{exchange}
Let $P$ be a nonempty finite set of integer  vectors in $\RR^n_+$ which contains with each 
$u\in P$ all its integral subvectors, and let $B(P)$ be the set of vectors $u\in P$ with $u<v$ 
for no $v\in P$. The following conditions are equivalent:
\begin{enumerate}
\item[(a)] $P$ is a discrete polymatroid;
\item[(b)] if $u,v\in P$ with $|v|>|u|$, then there is an integer $i$ such that 
$u+\epsilon_i\in P$ and 
\[
u+\epsilon_i\leq u\vee v;
\]
\item[(c)] 
\begin{enumerate}
\item[(i)] all $u\in B(P)$ have the same modulus,
\item[(ii)] if $u,v\in B(P)$ with $u(i)>v(i)$, then there exists $j$ with $u(j)<v(j)$ such 
that $u-\epsilon_i+\epsilon_j\in B(P)$.
\end{enumerate}
\end{enumerate}
\end{Theorem}

\begin{pf} That (a) implies (b) was already shown in the proof of \ref{cut}, and  implication 
(b)\implies (a) is trivial. 

(b)\implies (c): We noticed already that (c)(i) holds. Thus it remains to prove (c)(ii). 
Let $u,v\in B(P)$ with $u(i)>v(i)$ for some $i$. Then $u(i)-1\geq v(i)$, and hence 
$|u-\epsilon_i|=|v|-1<|v|$. Thus by (b) there exists an integer $j$ such that 
$(u-\epsilon_i)+\epsilon_j\leq (u-\epsilon_i)\wedge v$. We have $j\neq i$, because otherwise 
$u(i)=(u-\epsilon_i+\epsilon_j)(i)\leq \max\{u(i)-1,v(i)\}=u(i)-1$, a contradiction. 
Thus $u(j)+1=(u-\epsilon_i+\epsilon_j)(j)\leq \max\{u(j),v(j)\}\leq v(j)$, that is, 
$v(j)>u(j)$. 

(c)\implies (b): Let $v,u\in P$ with $|v|>|u|$, and let $w'\in B(P)$ with $u<w'$.  Since all 
$w'\in B(P)$ have the same modulus, it follows that $|v|\leq |w'|$. Thus we can choose a 
subvector $w$ of $w'$ in $P$ with $u\leq w$ and $|w|=|v|$. Let $P'=\{x\in P\: |x|\leq   
|v|\}$. By \ref{cut}(a), $P'$ is a discrete polymatroid, and hence by the implication 
(b)\implies (c) which we have already shown, we know that $w$ and $v$ satisfy the exchange 
property (c)(ii). 

Suppose that $w(j)\leq \max\{u(j),v(j)\}$ for all $j$. Since $u(j)\leq w(j)$ for all $j$, it 
follows that $w(j)\leq v(j)$ for all $j$. However since $|w|=|v|$, this implies $w=v$. Then 
$u<v$, and the assertion is trivial. 

Now assume that there exists an integer $j$ such that $w(j)>\max\{u(j),v(j)\}$. Then by the 
exchange property, there exists an integer $i$ with $v(i)>w(i)$ such that 
$w-\epsilon_i+\epsilon_j\in P$. Since $u+\epsilon_i$ is a subvector of 
$w-\epsilon_j+\epsilon_i$, it follows that $u+\epsilon_i\in P$. Furthermore we have 
$(u+\epsilon_i)(i)=u(i)+1\leq w(i)+1\leq v(i)$, so that $u+\epsilon_i\leq u\wedge v$.
\end{pf}

The following corollary is immediate from \ref{exchange} and the definition of a matroid.  

\begin{Corollary}
\label{matroid}
Let $B$ be a nonempty finite set of integer vectors in $\RR^n_{+}$. The following conditions 
are equivalent:
\begin{enumerate}
\item[(a)] $B$ is the set of bases of a matroid;
\item[(b)] $B$ is the set of bases of a discrete polymatroid, and for all $u\in B$ one has 
$u(i)\leq 1$ for 
$i=1,\ldots, n$.
\end{enumerate}
\end{Corollary}

The exchange property \ref{exchange}(c)(ii) suggests the following

\begin{Definition}
\label{defexchange}
{\em Let $B$ be a nonempty set of integer  vectors which have  the same modulus. Then $B$
satisfies: 
\begin{enumerate}
\item[(W)] the {\em weak exchange property}, if for all $u,v\in B$, $u\neq v$, there exist $i$ 
and $j$ with $u(i)>v(i)$ and $u(j)<v(j)$ such that $u-\epsilon_i+\epsilon_j\in B$,
\item[(S)] the {\em strong exchange property}, if for all $u,v\in B$, $u\neq v$,  and all $i$ 
and $j$ with $u(i)>v(i)$ and $u(j)<v(j)$ one has $u-\epsilon_i+\epsilon_j\in B$.
\end{enumerate}}
\end{Definition}

\begin{Examples}
\label{examples1}{\em 
(a) Let $B$ be the set of bases of a polymatroid on the ground set $[n]$ with $n\leq 3$. It is 
immediate that $B$ satisfies the strong exchange property.

(b) The  set $\{(1,1,1,1), (0,2,0,2), (0,1,1,2), (1,2,0,1)\}$ is the set of bases of a 
discrete polymatroid.  It does not satisfy the strong exchange property.

(c) Let $s_1,\ldots,s_n$ and $d$ be nonnegative integers. The set 
\[
V=\{u\: \text{$u(i)$ is an integer with $0\leq u(i)\leq s_i$ and $|u|=d$}\}
\]
is called of{\em Veronese type}. It satisfies the strong exchange property.

In fact, let $u,v\in V$ with $u(i)>v(i)$ and $u(j)<v(j)$ then $u(i)-1\geq 0$ and $u(j)+1\leq 
s_j$, so that $u-\epsilon_i+\epsilon_j\in V$.

(d) The set $\{(2,1,1), (2,2,0),(3,0,1),(3,1,0),(4,0,0)\}$ satisfies the strong exchange 
property, but is not of Veronese type.

(e) A set $S$ of integral vectors in $\RR^n_{+}$ of modulus $d$ is called {\em strongly 
stable}, if for all $u\in S$, all $i$ with $u(i)>0$ and all $j<i$, one has that 
$u-\epsilon_i+\epsilon_j\in S$. 
The strongly stable set $\{(3,0,1), (1,3,0), (3,1,0), (2,2,0),(4,0,0)\}$ does not  satisfies 
the weak exchange property. But any strongly stable set $S$ satisfies the following somewhat 
weaker exchange property: for all $u,v\in S$ there exist integers $i$ and $j$ such that  
$u(i)>v(i)$ and $u(j)<v(j)$,  and either $u-\epsilon_i+\epsilon_j\in S$, or
$v-\epsilon_j+\epsilon_i\in S$. 
 
Indeed, let $u,v\in S$, $u\neq v$. Then there is an integer $i$ such that $u(i)>v(i)$. Let $i$ 
be the largest such number. If $u(j)=v(j)$ for all $j>i$, then, since $u$ and $v$ have the 
same modulus, there exists $j<i$ with $u(j)<v(j)$. Since $S$ is strongly stable it follows 
that $u-\epsilon_i+\epsilon_j\in S$. On the other hand, if $u(j)<v(j)$ for some $j>i$, then 
$v-\epsilon_j+ \epsilon_i\in S$, since $S$ is strongly stable.

(f) A strongly stable set $S$ is called {\em principal}, if
there exists an element  $u\in S$ such that the smallest strongly stable set containing $u$ 
coincides with $S$. In that case $u$ is called the Borel generator of $S$. The set in example 
(d) is a  strongly stable principal set with Borel generator $(2,1,1)$. 

A strongly stable principal set satisfies the exchange property of \ref{exchange}(c)(ii), 
i.e.\ it is the set of bases of a polymatroid. This will follow from Example 
\ref{transexample} in Section 8. On the other hand, the set   
\[\{(0,1,0,1),(0,1,1,0),(0,2,0,0),(1,0,0,1),
(1,0,1,0),(1,1,0,0),(2,0,0,0)\}
\]
is strongly stable principal with    
Borel generator $(0,1,0,1)$, but does not satisfy the strong exchange property.   
}
\end{Examples}

We close this section with the following

\begin{Theorem}
\label{also}
Let $P$ be a nonempty finite set of integer vectors in $\RR^n_+$ which contains with each 
$u\in P$ all its integral subvectors. The following conditions are equivalent:
\begin{enumerate}
\item[(a)] $P$ is a discrete polymatroid of rank $d$ on the ground set $[n]$;
\item[(b)] $B=\{(u,d-|u|)\: u\in P\}$ is the set of bases  of a discrete polymatroid of rank 
$d$ on the ground set $[n+1]$.   
\end{enumerate}
\end{Theorem}

\begin{pf} (a)\implies(b): We will show that $B$ satisfies condition (c) of \ref{exchange}.  
Let $u,v\in P$, $i=d-|u|$,  $j=d-|v|$ and set $u'=(u,i)$ and $v'=(v,j)$. We may suppose that 
$|v|\geq |u|$, so that  $j\leq i$. If $i=j$, then $u$ and $v$ are bases of $P'=\{w\in P\: 
|w|\leq d-i\}$, see \ref{cut}. Thus $u'$ and $v'$ satisfy the exchange property  
\ref{exchange}(c)(ii), and hence we may  assume that $j<i$.
We consider two  cases. 

In the first case assume that $v'(k)>u'(k)$ for some $k$. Then $k\leq n$, and $v-\epsilon_k\in 
P$, since $v-\epsilon_k$ is a subvector of $v$, and so $v'-\epsilon_k+\epsilon_{n+1}\in B$. 
  
In the second case assume that $u'(k)>v'(k)$ for some $k$. Since $|v|>|u|$, \ref{exchange}(b) 
implies that there exists an integer $l$ such that $u+\epsilon_l\in P$ and $u+\epsilon_l\leq 
u\wedge v$. It follows that $u(l)<u(l)+1\leq v(l)$. If $k\leq n$, then  
$u-\epsilon_k+\epsilon_l\in P$, because it is a subvector of $u-\epsilon_l$. Hence  
$u'-\epsilon_k+\epsilon_l=(u-\epsilon_k+\epsilon_l, i)\in B$. 
On the other hand, if $k=n+1$, that is, $u'(k)=i$, then 
$u'-\epsilon_k+\epsilon_l=(u+\epsilon_l,i-1)\in B$, because $u+\epsilon_l\in P$.

(b)\implies (a): Let $u,v\in P$ with $|v|>|u|$. Then $d-|v|<d-|u|$, and so by the exhcange 
property of $B$ there exists an integer $i$ with $v(i)>u(i)$ and such that 
$(u+\epsilon_i,d-|u|)\in B$. This implies that $u+\epsilon_i\in P$, and since $v(i)>u(i)$ we 
also have that $u+\epsilon_i\leq u\wedge v$.
\end{pf}

\section{Integral polymatroids and discrete polymatroids}
The purpose of the present section is to show our first
fundamental Theorem \ref{firstmaintheorem}
which says that a nonempty finite set
$P \subset \ZZ_{+}^{n}$ is a discrete polymatroid
if and only if $\conv(P) \subset \RR_{+}^{n}$
is an integral polymatroid with
$\conv(P) \cap \ZZ^{n} = P$.
Here $\conv(P)$ is the convex hull of $P$ in $\RR^n$.

First of all, we collect a few basic lemmata on integral
polymatroids and discrete polymatroids which will be
required to prove Theorem \ref{firstmaintheorem}.

\begin{Lemma}
\label{integralpolymatroid}
If ${\cal P} \subset \RR_{+}^{n}$ is an integral
polymatroid and if $u, v \in {\cal P} \cap \ZZ^n$
with $|v| > |u|$, then there is
$w \in {\cal P} \cap \ZZ^n$ such that
$u < w \leq u \vee v$.
\end{Lemma}

Lemma \ref{integralpolymatroid}
appears (not explicitly) in
\cite[Lemma 5, p.\ 340]{Welsh}.
Its proof is also valid to show
Lemma \ref{integralpolymatroid}
by noting that the ground set rank function $\rho$
of an integral polymatroid ${\cal P}$
is integer valued.

Let $P \subset \ZZ_+^n$ be a discrete polymatroid
and $B(P)$ the set of bases of $P$.
We define the nondecreasing function
$\rho_P \: 2^{[n]} \to \RR_+$
associated with $P$ by setting
\[
\rho_P(X) = \max\{ u(X) \: u \in B(P) \}
\]
for all $\emptyset \neq X \subset [n]$
together with $\rho_P(\emptyset) = 0$.

In general, for $u, v \in B(P)$ we define the
{\em distance} between $u$ and $v$ by
\[
\dis(u,v) =
\frac{1}{2}
\sum_{i=1}^{n} |u(i) - v(i)|.
\]
A crucial property of $\dis(u,v)$ is that
if $u(i) > v(i)$ and $u(j) < v(j)$
together with
$u' = u - \varepsilon_i + \varepsilon_j \in B(P)$,
then $\dis(u,v) > \dis(u',v)$.

\begin{Lemma}
\label{discretepolymatroid}
If $X_1 \subset X_2 \subset \cdots \subset X_s
\subset [n]$ is a sequence of subsets of $[n]$,
then there is $u \in B(P)$ such that
$u(X_k) = \rho_P(X_k)$ for all $1 \leq k \leq s$.
\end{Lemma}

\begin{pf}
We work with induction on $s$ and suppose that
there is $u \in B(P)$ such that
$u(X_k) = \rho_P(X_k)$ for all $1 \leq k < s$.
Choose $v \in B(P)$ with $v(X_s) = \rho_P(X_s)$.
If $u(X_s) < v(X_s)$, then there is
$i \in [n]$ with $i \not\in X_s$
such that $u(\{i\}) > v(\{i\})$.
The exchange property \ref{exchange} (c) (ii)
says that there is
$j \in [n]$ with $u(j) < v(j)$
such that
$u_1
= u - \varepsilon_{i} + \varepsilon_{j}
\in B(P)$.  Since $u(X_{s-1}) = \rho_P(X_{s-1})$,
it follows that $j \not\in X_{s-1}$.
Hence $u_1(X_k) = \rho_P(X_k)$ for all
$1 \leq k < s$.  Moreover,
$u_1(X_s) \geq u(X_s)$
and $\dis(u,v) > \dis(u_1,v)$.

If $u_1(X_s) = v(X_s)$, then $u_1$ is a desired
base of $P$.  If $u_1(X_s) < v(X_s)$, then
the above technique will yield a base $u_2$ of
$P$ such that
$u_2(X_k) = \rho_P(X_k)$ for all
$1 \leq k < s$,
$u_2(X_s) \geq u_1(X_s)$
and $\dis(u_1,v) > \dis(u_2,v)$.
It is now clear that repeated applications of
this argument guarantee the existence of
a base $u_q$ of $P$ such that
$u_q(X_k) = \rho_P(X_k)$ for all $1 \leq k \leq s$.
\end{pf}

\begin{Corollary}
\label{submodularassociatedwith}
The function
$\rho_P \: 2^{[n]} \to \RR_+$
is submodular.
\end{Corollary}

\begin{pf}
Let $A, B \subset [n]$.
By Lemma \ref{discretepolymatroid}
there is $u \in B(P)$ such that
$u(A \cap B) =\rho_P(A \cap B)$
and
$u(A \cup B) =\rho_P(A \cup B)$.
Hence
\begin{eqnarray*}
\rho_P(A) + \rho_P(B) & \geq & u(A) + u(B) \\
& = & u(A \cup B) + u(A \cap B) \\
& = & \rho_P(A \cup B) + \rho_P(A \cap B),
\end{eqnarray*}
as desired.
\end{pf}

We now come to our first fundamental

\begin{Theorem}
\label{firstmaintheorem}
A nonempty finite set
$P \subset \ZZ_{+}^{n}$ is a discrete polymatroid
if and only if $\conv(P) \subset \RR_{+}^{n}$
is an integral polymatroid with
$\conv(P) \cap \ZZ^{n} = P$.
\end{Theorem}

\begin{pf}
The ``\,if\,'' part follows from
Lemma \ref{integralpolymatroid}.
To see why the ``\,only if\,'' part is true, let
$P \subset \ZZ_{+}^{n}$ be a discrete polymatroid
and $\rho_P \: 2^{[n]} \to \RR_+$
the nondecreasing and submodular function
associated with $P$.  Write
${\cal P} \subset \RR_{+}^{n}$ for the integral
polymatroid with $\rho_P$ its ground set rank
function, i.e.,
\[
{\cal P} = \{ u \in \RR_+^n \:
u(X) \leq \rho_P(X), X \subset [n] \}.
\]
Since each base $u$ of $P$ satisfies
$u(X) \leq \rho_P(X)$ for all $X \subset [n]$,
it follows that $P \subset {\cal P}$.  Moreover,
since ${\cal P}$ is convex, one has
$\conv(P) \subset {\cal P}$.
Now, Lemma \ref{discretepolymatroid}
together with Proposition \ref{vertex}
guarantees that all vertices of
${\cal P}$ belong to $P$.  Thus
${\cal P} = \conv(P)$.

To complete our proof we must show
${\cal P} \cap \ZZ^n = P$.
For each $i \in [n]$,
write $P_i \subset \ZZ_+^n$
for the discrete polymatroid
$P_{\varepsilon_i}$
in the notation of Lemma \ref{cut} (b)
and $B_i = B(P_i)$, the set of bases of
$P_i$.
We compute
the nondecreasing and submodular function
$\rho_{P_i} \: 2^{[n]} \to \RR_+$
associated with $P_i$.
We distinguish three cases:

\medskip
(a)
Let $i \not\in X \subset [n]$ with
$\rho_P(X \cup \{i\}) > \rho_P(X)$.
Choose
$u \in B(P)$ with $u(X) = \rho_P(X)$ and with
$u(X \cup \{i\}) = \rho_P(X \cup \{i\})$.
Then $u(i) = u(X \cup \{i\}) - u(X) \geq 1$
and $u \in B(P)$.  Since $i \not\in X$, one has
$(u - \varepsilon_i)(X) = u(X)$.
Since $(u - \varepsilon_i)(X) \leq
\rho_{P_i}(X) \leq \rho_P(X) = u(X)$,
it follows that
$\rho_{P_i}(X) = \rho_P(X)$.

\medskip
(b)
Let $i \not\in X \subset [n]$ with
$\rho_P(X \cup \{i\}) = \rho_P(X)$.
If $u \in B(P)$ with $u(i) \geq 1$, then
\begin{eqnarray*}
(u - \varepsilon_i)(X) & \leq &
(u - \varepsilon_i)(X \cup \{i\}) \\
& = & u(X \cup \{i\}) - 1 \\
& \leq &
\rho_P(X \cup \{i\}) - 1
= \rho_P(X) - 1.
\end{eqnarray*}
Thus
$\rho_{P_i}(X) \leq \rho_P(X) - 1$.
Choose $v \in B(P)$ with $v(X) = \rho_P(X)$.
Then $v(i) = 0$.
(Otherwise, since $i \not\in X$,
$\rho_P(X \cup \{i\}) \geq v(X \cup \{i\})
= v(X) + v(i) > v(X) = \rho_P(X)$,
a contradiction.)
Let $u_0 \in B(P)$ with $u_0(i) \geq 1$.
Then the exchange property says that
there is $j \in [n]$ with
$u_0(j) < v(j)$ such that
$u_0 - \varepsilon_i + \varepsilon_j \in B(P)$.
Thus we assume $u_0(i) = 1$.
If $u_0(X) < v(X) - 1$, then
$u_0(X \cup \{i\}) = u_0(X) + 1 <
v(X) = v(X \cup \{i\})$.  Thus there is
$j \not\in X \cup \{i\}$ with
$u_0(j) > v(j)$.  Hence there is
$i \neq k \in [n]$ with
$u_0(k) < v(k)$
such that
$u_1 = u_0 - \varepsilon_j + \varepsilon_k
\in B(P)$.  Then $u_1(i) = 1$,
$u_1(X) \geq u_0(X)$ and
$\dis(u_0,v) > \dis(u_1,v)$.
Thus, as in the proof of
Lemma \ref{discretepolymatroid},
we can find $u \in B(P)$ with
$u(i) = 1$ such that
$u(X) = (u - \varepsilon_i) (X) = v(X) - 1$.
Hence
$\rho_{P_i}(X) = \rho_P(X) - 1$.

\medskip
(c)
Let $i \in X \subset [n]$.
Then
$\rho_{P_i}(X) = \rho_P(X) - 1$.
In fact, since $i\in X$,
by Lemma \ref{discretepolymatroid}
there is $u \in B(P)$ with
$u(i) = \rho_P(\{i\}) \geq 1$
and with $u(X) = \rho_P(X)$.
Then $(u - \varepsilon_i)(X)
= \rho_P(X) - 1$.

\medskip

Let ${\cal P}_i \subset \RR_+^n$
denote the integral polymatroid
with $\rho_{P_i}$ its ground set rank
function.  Then
${\cal P}_i = \conv(P_i)$ and,
working with induction on the rank of $P$
enables us to assume that
${\cal P}_i \cap \ZZ^n = P_i$.
If $x \in {\cal P} \cap \ZZ^n$ with
$x(i) \geq 1$, then $y = x - \varepsilon_i$
belongs to ${\cal P}_i$.  (In fact,
if $i \not\in X \subset [n]$ with
$\rho_{P_i}(X) = \rho_P(X) - 1$,
then
$\rho_P(X \cup \{i\}) = \rho_P(X)$.
Thus replacing $u$ with $x$ in the inequalities
$(u - \varepsilon_i)(X) \leq
\cdots = \rho_P(X) - 1$
appearing in the discussion (b)
shows that
$y(X) \leq \rho_{P_i}(X)$.)
Thus $y \in {\cal P}_i \cap \ZZ^n = P_i$.
Hence $y \leq u - \varepsilon_i$
for some $u \in B(P)$ with
$u(i) \geq 1$.  Thus $x \leq u \in B(P)$
and $x \in P$.  Hence
${\cal P} \cap \ZZ^n = P$,
as desired.
\end{pf}
\section{The symmetric exchange theorem for polymatroids}
We now establish our second fundamental theorem
on discrete polymatroids;
the symmetric exchange theorem.

\begin{Theorem}
\label{secondmaintheorem}
If $u = (a_1, \ldots, a_n)$ and
$v = (b_1, \ldots, b_n)$ are bases of a discrete
polymatroid $P \subset \ZZ_+^n$, then
for each $i \in [n]$ with $a_i > b_i$
there is $j \in [n]$ with $a_j < b_j$
such that both $u - \varepsilon_i + \varepsilon_j$
and $v - \varepsilon_j + \varepsilon_i$
are bases of $P$.
\end{Theorem}

\begin{pf}
Let $B'$ denote the set of those bases $w$ of $P$
with $u \wedge v \leq w \leq u \vee v$.
It then turns out that $B'$ satisfies the exchange
property \ref{exchange} (c) (ii) for polymatroids.
Thus $B'$ is the set of bases
of a discrete polymatroid $P' \subset \ZZ_+^n$.
Considering $u' = u - u \wedge v$ and
$v' = v - u \wedge v$ instead of $u$ and $v$,
we will assume that $P' \subset \ZZ_+^s$
is a discrete polymatroid, where $s \leq n$, and
\[
u = (a_1, \ldots, a_r, 0, \ldots, 0) \in \ZZ_+^s,
\, \, \, \, \,
v = (0, \ldots, 0, b_{r+1}, \ldots, b_s)
\in \ZZ_+^s,
\]
where each $0 < a_i$ and each $0 < b_j$
and where $|u| = |v| = \rank(P')$.
Our work is to show that
for each $1 \leq i \leq r$
there is $r + 1 \leq j \leq s$
such that
both $u - \varepsilon_i + \varepsilon_j$
and $v - \varepsilon_j + \varepsilon_i$
are bases of $P'$.
Let, say, $i = 1$.

\medskip
\noindent
First case: Suppose that
$u - \varepsilon_1 + \varepsilon_j$
are bases of $P'$ for all $r + 1 \leq j \leq s$.
It follows from the exchange property that,
given arbitrary $r$ integers $a'_1, \ldots, a'_r$
with each $0 \leq a'_i \leq a_i$,
there is a base $w'$ of $P'$ of the form
\[
w' = (a'_1, \ldots, a'_r, b'_{r+1}, \ldots, b'_s),
\]
where each $b'_j \in \ZZ$
with $0 \leq b'_j \leq b_j$.
Thus in particular there is
$r + 1 \leq j_0 \leq s$ such that
$v - \varepsilon_{j_0} + \varepsilon_1$
is a base of $P'$.  Since
$u - \varepsilon_1 + \varepsilon_j$
is a base of $P'$ for each
$r + 1 \leq j \leq s$,
both
$u - \varepsilon_1 + \varepsilon_{j_0}$
and
$v - \varepsilon_{j_0} + \varepsilon_1$
are bases of $P'$, as desired.

\medskip
\noindent
Second case: Let $r \geq 2$ and $r + 2 \leq s$.
Suppose that there is $r + 1 \leq j \leq s$
with
$u - \varepsilon_1 + \varepsilon_j
\not\in P'$.
Let $X \subset \{r + 1, \ldots, s\}$
denote the set of those
$r + 1 \leq j \leq s$
with
$u - \varepsilon_1 + \varepsilon_j
\not\in P'$.
Recall that Theorem \ref{firstmaintheorem}
guarantees that $\conv(P') \subset \RR_+^s$
is an integral polymatroid on the ground set
$[s]$ with $\conv(P') \cap \ZZ^s = P'$.
Let $\rho = \rho_{P'}$
denote the ground set rank function
of the integral polymatroid
$\conv(P') \subset \RR_+^s$.
Thus
$\rho(Y) = \max\{w(Y) \: w \in B' \}$
for $\emptyset \neq Y \subset [s]$
together with $\rho(\emptyset) = 0$.
In particular
$\rho(Y) = u(Y)$ if
$Y \subset \{1, \ldots, r\}$ and
$\rho(Y) = v(Y)$ if
$Y \subset \{r + 1, \ldots, s\}$.

For each $j \in X$, since
$u - \varepsilon_1 + \varepsilon_j
\not\in \conv(P')$, there is a subset
$A_j \subset \{2, 3, \ldots, r\}$
with
\[
\rho(A_j \cup \{j\}) \leq u(A_j).
\]
Thus
\begin{eqnarray*}
\rho(\{2, 3, \ldots, r\} \cup \{j\})
& \leq &
\rho(A_j \cup \{j\}) +
\rho(\{2, 3, \ldots, r\} \setminus A_j) \\
& \leq & u(A_j) +
u(\{2, 3, \ldots, r\} \setminus A_j) \\
& = & u(\{2, 3, \ldots, r\}) =
\rho(\{2, 3, \ldots, r\}).
\end{eqnarray*}
Hence, for all $j \in X$,
\begin{eqnarray*}
\rho(\{2, 3, \ldots, r\} \cup \{j\})
= u(\{2, 3, \ldots, r\}).
\end{eqnarray*}
By \cite[Lemma 1.3.3]{Oxley}
it follows that
\begin{eqnarray*}
\rho(\{2, 3, \ldots, r\} \cup X)
= u(\{2, 3, \ldots, r\}).
\end{eqnarray*}
Now, since $\rho$ is submodular,
\begin{eqnarray*}
\rho(\{2, 3, \ldots, r\} \cup X) +
\rho(\{1\} \cup X)
& \geq &
\rho(X) +
\rho(\{1, 2, \ldots, r\} \cup X) \\
& = & v(X) + \rank(P').
\end{eqnarray*}
Thus
\[
u(\{2, 3, \ldots, r\}) + \rho(\{1\} \cup X)
\geq v(X) + \rank(P').
\]
Since
\begin{eqnarray*}
\rank(P') - u(\{2, 3, \ldots, r\}) = a_1,\quad \text{and}
\\
\rho(\{1\} \cup X)
\leq \rho(\{1\}) + \rho(X)
= a_1 + v(X),
\end{eqnarray*}
it follows that
\[
\rho(\{1\} \cup X)
= a_1 + v(X).
\]
Hence, for all $X' \subset X$,
\begin{eqnarray*}
a_1 + v(X)
& = &
a_1 + v(X') + v(X \setminus X') \\
& = & \rho(\{1\}) + \rho(X')
+ \rho(X \setminus X') \\
& \geq & \rho(\{1\} \cup X') + \rho(X \setminus X') \\
& \geq & \rho(\{1\} \cup X) \\
& = & a_1 + v(X).
\end{eqnarray*}
Thus, for all $X' \subset X$,
\[
\rho(\{1\} \cup X')
= a_1 + v(X').
\]
By virtue of Lemma \ref{discretepolymatroid}
there is a base $w$ of $P'$
with $w(1) = a_1$ and with
$w(j) = v(j)$
$( = \rho(\{j\}) )$
for all $j \in X$.
Again the exchange property (for $w$ and $v$)
guarantees that for each $1 \leq i \leq r$
with $w(i) > 0$ there is
$j \in \{r + 1, \ldots, s\} \setminus X$
such that
$w - \varepsilon_i + \varepsilon_j$
is a base of $P'$.  Hence
repeated applications of the exchange property
yield a base of $w'$ of $P'$
of the form
$w' = v - \varepsilon_{j_0} + \varepsilon_1$,
where
$j_0 \in \{r + 1, \ldots, s\} \setminus X$.
Hence both
$u - \varepsilon_1 + \varepsilon_{j_0}$
and
$v - \varepsilon_{j_0} + \varepsilon_1$
are bases of $P'$, as required.
\end{pf}

\section{Base rings of  discrete polymatroids and their relations}
 
Let $K$ be a field, and $P$  a discrete polymatroid with the set of bases $B$. The toric ring 
$K[B]$  generated over $K$ by the monomials $t^u=\prod_it_i^{u(i)}$, $u=(u(1),\ldots,u(n))\in 
B$,  is called the {\em base ring of $P$}.  Let $S=K[x_u\: u\in B]$ be the polynomial ring in 
the indeterminates $x_u$ with $u\in B$. Let  $I_B$ be the kernel of the $K$-algebra 
homomorpism $\xi\: S\to K[B]$ with $\xi(x_u)=t^u$. There are some obvious elements in $I_B$, 
namely, those arising from symmetric exhange: Let $u, v\in B$ with $u(i)>v(i)$ and 
$u(j)<v(j)$, and such that $u-\epsilon_i+\epsilon_j$ and $v-\epsilon_j+\epsilon_i$ belong to 
$B$. Then clearly $x_ux_v-x_{u-\epsilon_i+\epsilon_j}x_{v-\epsilon_j+\epsilon_i}\in I_B$. We 
call such a relation a {\em symmetric exchange relation}. White conjectured (cf.\ 
\cite{White}) that for a matroid the symmetric exchange relations generate $I_B$. It is 
natural to conjecture that this also holds  for discrete polymatroids. 

In order to formulate the next result we have to recall the notion of sortability, introduced 
by Sturmfels \cite{Sturmfels}:

Let $u,v\in B$, and write $t^ut^v=t_{i_1}\ldots t_{i_{2d}}$ with $i_1\leq i_2\leq \cdots \leq 
i_{2d}$. Then we set $t^{u'}=\prod_j^dt_{2j-1}$ and $t^{v'}=\prod_j^dt_{2j}$. This defines a 
map 
\[
\sort\: B\times B\to M_d\times M_d,\quad (u,v)\to (u',v'),
\]
where $M_d$ denotes the set of all integer vectors of modulus $d$. 

The map `$\sort$' is called  the {\em sorting operator}, and $B$ is called {\em sortable},
if $\Im(\sort)\subseteq B\times B$. 

The pair $(u,v)$ is called {\em sorted} if $(u,v)\in \Im(\sort)$, equivalently, if 
$\sort(u,v)=(u,v)$. It is clear from the definition that $(u,v)$ is sorted if and only if 
\[
u(1)+\cdots +u(i)-1\leq v(1)+\ldots +v(i)\leq u(1)+\ldots +u(i)\quad\text{for $i=1,\ldots n$.} 
\]
This implies that $u-v$ is vector with entries $\pm 1$ and $0$. To such a vector we attach 
sequence of $+$ and $-$ signs: reading from the left to right we put the sign $+$ or the sign 
$-$ if we reach the entry $+1$ or the entry $-1$. For example to $(0,1,1,0,-1,-1)$ belongs the 
sequence $++--$.

For the above characterization of sorted pairs one easily deduces 

\begin{Lemma}
\label{sorted}
The following conditions are equivalent:
\begin{enumerate}
\item[(a)] $(u,v)$ is sorted;
\item[(b)] $u-v$ is a vector with entries $\pm 1$ and $0$ and with $\pm$-sequence 
$+-+-\cdots+-$.
\end{enumerate}
\end{Lemma}

In \cite{DeNegri}, De Negri noticed the following fact which easily follows from a theorem of 
Sturmfels \cite{Sturmfels}.

\begin{Lemma}
\label{restriction}
Suppose $B\subset M_d$ is sortable. Then $I_B$ has a Gr\"obner basis consisting 
of the sorting relations $x_ux_v-x_{u'}x_{v'}$ with $u,v\in B$ and $(u',v')=\sort(u,v)$.  
\end{Lemma}

The main result of this section is

\begin{Theorem}
\label{partial} 
{\em (a)} Suppose that each matroid has the property that the toric ideal of its
base ring is generated by  symmetric exchange relations, then this is also true for each 
discrete polymatroid. 

{\em (b)} If $P$ is a discrete polymatroid whose set of bases $B$ satisfies the strong 
exchange property, then 
\begin{enumerate}
\item[(i)] $B$ is sortable, and hence $I_B$ has a quadratic Groebner basis and $K[B]$ is 
Koszul,
\item[(ii)] $I_B$ is generated by  symmetric exchange relations. 
\end{enumerate}
\end{Theorem}

For the proof of \ref{partial}(a) we shall need

\begin{Lemma}
\label{exchangerel}
Let $u_1,\ldots,u_d\in B$. Then, modulo exchange relations, $\prod_j^dx_{u_j}$ can be 
rewritten as $\prod_j^dx_{v_j}$ with $|v_j(i)-v_k(i)|\leq 1$ for all $i,j, k$.
\end{Lemma}

\begin{pf}
Suppose that for some $i$, $j$ and $k$ we have $|u_k(i)-u_l(i)|>1$. Without loss of generality 
we may assume that $k=1$, $l=2$ and $u_1(i)>u_2(i)$. By the symmetric exchange property there 
exists $j$ with $u_2(j)>u_1(j)$ such that $u_1-\epsilon_i+\epsilon_j, 
u_2-\epsilon_j+\epsilon_i\in B$.

We set 
\[
u_k^{'}=\left\{
\begin{array}{ll}
u_k, &\text{for $k\neq 1,2$},\\
u_1-\epsilon_i+\epsilon_j, & \text{for $k=1$},\\
u_2-\epsilon_j+\epsilon_i, & \text{for $k=1$}.\\
\end{array}
\right. 
\]
It is clear that $u_1^{'}(i)-u'_2(i)=u_1(i)-u_2(i)-2\geq 0$.

We introduce the number
\[
c_i(u_1,\ldots,u_d)=\sum_{1\leq k<l\leq d}|u_k(i)-u_l(i)|.
\]
Then $c_k(u_1,\ldots,u_d)=c_k(u'_1,\ldots,u'_d)\quad\text{for $k\neq i,j$}.$

It is easy to show that 
\begin{eqnarray}
\label{three} 
c_i(u'_1,\ldots,u'_d)<c_i(u_1,\ldots,u_d),
\end{eqnarray}
and
\begin{eqnarray}
\label{four}
c_j(u'_1,\ldots,u'_d)\leq c_j(u_1,\ldots,u_d).
\end{eqnarray}
Thus, induction completes the proof.
\end{pf}

\begin{pf}[ Proof of \ref{partial}] 
(a) Let $\prod_jx_{u_j}-\prod_jx_{v_j}\in I_B$. Applying Lemma \ref{exchangerel} we may assume 
that for all $j$ and $k$, all components of $u_j$ and $u_k$, and of $v_j$ and $v_k$ differ at 
most by $1$.  Let $w=\sum_j^du_j=\sum_j^dv_j$, and for $i=1,\ldots n$ let $a_i=\min\{u_j(i)\: 
j=1,\ldots, d\}$, and similarly $b_i= \min\{v_j(i)\: j=1,\ldots, d\}$.
Then $u_j(i)=a_i$ or $u_j(i)=a_i+1$. Let $k_i$ be the number of $j$ with $u_j(i)=a_i+1$. Then  
$w(i)=\sum_{j=1}^du_j(i)= da_i+k_i$. Since $0\leq k_i<d$, we conclude that $a_i=\lfloor 
w(i)/d\rfloor$, where  $\lfloor c\rfloor$ denotes the integer part of a number $c$. In 
particular it follows that $a_i=b_i$ for $i=1,\ldots, n$. 

Now we consider the following subset 
\[
B'=\{u\in B\: a_i\leq u(i)\leq a_i+1\}
\]
of $B$. The vectors $u_1,\ldots,u_d$ and $v_1,\ldots, v_d$ belong to $B'$, and $B'$ is closed 
under the exchange relations. 

We may identify $B'$ with the set of bases of a matroid $B''$, via the assignment $u\mapsto 
u^*=u-(a_1,\ldots,a_n)$. Then the relation $\prod_jx_{u_j}-\prod_jx_{v_j}\in I_B$ corresponds 
to the relation $\prod_jx_{u^*_j}-\prod_jx_{v^*_j}\in I_{B''}$. By our hypothesis, 
$\prod_jx_{u^*_j}-\prod_jx_{v^*_j}$ reduces to $0$ modulo the  symmetric
exchange relations of $B''$. Thus $\prod_jx_{u_j}-\prod_jx_{v_j}$ reduces to $0$ modulo the 
symmetric exchange relations of $B'$, and hence of $B$.

(b) (i) Let $(u,v)\in B\times B$. In the proof of (a) we have shown that 
by a finite number of exchanges we can transform the pair $(u,v)$ into a pair whose difference 
vector has entries $\pm 1$ and $0$. Thus we may assume that $(u,v)$ itself has this property. 
Since $u$ and $v$ have the same modulus, it follows that the $\pm 1$-sequence of $u-v$ has as 
many $+$ signs as $-$ signs. If the $\pm$-sequence does not have the pattern $+-+-\cdots+-$, 
then obviously we can  obtain such a sequence by a finite number of symmetric exchanges. Thus 
\ref{sorted} shows that $B$ is sortable. The other statements of (i) follow from 
\ref{restriction}.

(ii) We have seen in the proof of (i) that by a finite number of symmetric exchanges we can 
sort any pair $(u,v)\in B$. It follows that the sorting relation $x_ux_v-x_{u'}x_{v'}$ is  a 
linear combination of the corresponding symmetric exchange relations. Since by (i), the 
sorting relations generate $I_B$, the assertion follows.  
\end{pf}

\section{The Ehrhart ring and the base ring of a discrete polymatroid}

Let $K$ be a field, and let $P$ be a discrete polymatroid of rank $d$ on the ground set $[n]$ 
with set of bases 
$B$. Since $P$ is the set of integer  vectors of an integral polymatroid ${\cal P}$ (see 
Theorem \ref{firstmaintheorem}), we may study the {\em Ehrhart ring of ${\cal P}$}. To this 
end one considers the cone ${\cal C}\subset \RR^{n+1}$ with ${\cal C}=\RR_+\{(p,1)\: p\in 
{\cal P}\}$. Then $Q={\cal C}\sect \ZZ^{n+1}$ is subsemigroup of $\ZZ^{n+1}$, and the Ehrhart 
ring of ${\cal P}$  is defined to be the  toric ring $K[{\cal P}]\subset K[t_1,\ldots,t_n,s]$ 
generated over $K$ by the monomials $t^us^i$, $(u,i)\in Q$. By Gordan's lemma (cf.\ 
\cite[Proposition 6.1.2]{BruHer}), $K[{\cal P}]$ is normal. 

Notice that $K[{\cal P}]$ is naturally graded if we assign to $t^us_i$ the degree $i$. We 
denote by  $K[P]$ the $K$-subalgebra of $K[{\cal P}]$ which is generated over $K$ by the 
elements  of degree 1 in
$K[{\cal P}]$. Since $P={\cal P}\sect \ZZ^n$  it follows that 
$K[P]=K[t^us\: u\in P]$. Observe that $K[B]$ may be identified with the subalgebra $K[t^us\: 
u\in B]$ of $K[P]$. After this identification $K[B]$ is an algebra retract of $K[P]$ with 
retraction map $\pi\: K[P]\to K[B]$, where $\pi$ is the residue class map modulo the prime 
ideal $\wp$ generated  by the elements $t^us$, $u\in P$, $|u|<d$.

As an immediate consequence of \ref{sum} we obtain the following important

\begin{Theorem}
\label{conclusion} 
$K[P]=K[{\cal P}]$. In particular, $K[P]$ is normal.
\end{Theorem}

\begin{Corollary}
\label{alsonormal}
$K[B]$ is normal. 
\end{Corollary}

\begin{pf}

By \cite[Theorem 6.1.4]{BruHer}, $K[P]$ is normal if and only if the semigroup $S$ generated 
by $\hat{P}=\{(u,1)\: u\in P\}$ is normal. This means that if $G$ is the smallest subgroup of 
$\ZZ^{n+1}$ containing $S$ and if $ma\in S$ for some $m\in \NN$ and $a\in G$, then $a\in S$. 
We apply the same criterion to show that $K[B]$ is normal.  
Let $H$ be the subgroup  of $\ZZ^{n+1}$ consisting of all $(u,i)$ with $|u|=id$, and let $T$ 
be the semigroup generated by  $\hat{B}=\{(u,1)\: u\in B\}$. Then $H\sect G$ is the smallest 
group containing the semigroup $T$. Hence if $ma\in T$ for some $m\in \NN$ and $a\in H\sect 
G$, then  $a\in S\sect H=T$, as desired.
\end{pf}

Our next aim is to compare algebraic properties of $K[P]$ with those of $K[B]$.
We say that a $K$-algebra $A$ has a quadratic Gr\"obner basis if the defining ideal of $A$ has 
this property for some term order.

\begin{Theorem}
\label{compare}
(a) Suppose $K[P]$ has quadratic relations, or a quadratic Gr\"obner basis or is Koszul, then 
$K[B]$ has these properties, too.

(b) Given a property ${\cal E}$. Suppose that  $K[B(P)]$ satisfies ${\cal E}$ for all discrete 
polymatroids $P$. Then also $K[P]$ satisfies ${\cal E}$ for all discrete polymatroids $P$.
\end{Theorem}

In particular it follows from the theorem that the properties that the Ehrhart ring of all 
discrete polytopes is Koszul, if and only if this is the case for all  base rings of all 
discrete polytopes. Not all properties behave this way. For example, there exist discrete 
polytopes (see \ref{GorEx}) for which $K[P]$ is Gorenstein but $K[B]$ is not, and vice versa.

\begin{pf}[Proof of \ref{compare}] (a) We noticed already  that $K[B]$ may be viewed as  an 
algebra retract of $K[P]$. Hence the statements concerning quadratic relations and Koszulness 
follow from \cite[Proposition 1.4, Corollary 2.5]{OhHerHi}.   

Now let $A\subset B$ be an arbitrary  retract of standard graded toric $K$-algebras. Then we 
may assume that $A=K[x_1,\ldots,x_n]/I$, and $B=K[x_1,\ldots,x_n,y_1,\ldots,y_m]/J$ with toric 
ideals $I$ and $J$ contained in the square of the corresponding graded maximal ideals, and 
that the retraction map $B\to A$ is induced by $\pi(y_j)=0$ for $j=1,\ldots,m$. 

We use the following well-known fact: suppose $J$ has a Gr\"obner basis $G$ with respect to 
the term order $<$, and suppose that for all $f\in G$ with $\inii_<(f)\in K[x_1,\ldots,x_n]$ 
one has that $f\in k[x_1,\ldots, x_n]$. Then $G'=G\sect K[x_1,\ldots,x_n]$ is a Gr\"obner 
basis of $I$ (with respect to the restricted term order).   

Since $B$ is toric, we may assume that the elements of  $G$ are binomials. Suppose now that 
$f\in G$ with $f=m_1-m_2$ and $\inii(f)=m_1\in K[x_1,\ldots,x_n]$, and suppose that 
$m_2\not\in 
K[x_1,\ldots,x_n]$. Then $m_1=\pi(f)\in I$, a contradiction since $I$ is a prime ideal.

(b) We notice that $K[P]$ is isomorphic to the toric ring $K[t^us^{d-|u|}\: u \in P]$ which, 
as we know from \ref{also}, is the base ring of a polymatroid. Thus the conclusion follows. 
\end{pf}

\begin{Remark}
{\em A monomial ideal generated by monomials corresponding to the base of a discrete 
polymatroid is called a {\em polymatroidal ideal}. It has been shown in \cite{HerTak} that 
polymatroidal ideals have linear quotients, so that, in particular, they have a linear 
resolution. As another immediate application of \ref{sum} one has that the product of 
polymatroidal ideals is again polymatroidal. This fact has been shown directly in 
\cite{ConHer}. 
}
\end{Remark}

\section{Gorenstein polymatroids}
Let $P \subset \ZZ_+^n$ be a discrete polymatroid
and $B = B(P)$ the set of bases of $P$.
We know that both the Ehrhart ring $K[P]$ and
the base ring $K[B]$ are normal; thus
Cohen--Macaulay.  It is then natural to ask
when these rings are Gorenstein.
We begin with

\begin{Example}
\label{GorEx}
{\em
(a) Let $P \subset \ZZ_+^3$ be the discrete
polymatroid consisting of all integer vectors
$u \in \ZZ_+^3$ with $|u| \leq 3$.  Then
the base ring $K[B]$ is the Veronese subring
$K[x, y, z]^{(3)}$, the subring of $K[x, y, z]$
generated by all monomials of degree $3$.  Thus
$K[B]$ is Gorenstein.  On the other hand, since
the Hilbert series of the Ehrhart ring $K[P]$
is $(1 + 16 \lambda + 10 \lambda^2)
/ (1 - \lambda)^4$, it follows that
$K[P]$ is not Gorenstein.

(b) Let $P \subset \ZZ_+^3$ be the discrete
polymatroid consisting of all integer vectors
$u \in \ZZ_+^3$ with $|u| \leq 4$.  Then
$K[B] = K[x, y, z]^{(4)}$ is not Gorenstein.
On the other hand,
the Hilbert series of the Ehrhart ring $K[P]$
is $(1 + 31 \lambda + 31 \lambda^2 + \lambda^3)
/ (1 - \lambda)^4$; thus
$K[P]$ is Gorenstein.

(c) Let $P \subset \ZZ_+^2$ be the discrete
polymatroid with $B = \{ (1, 2), (2, 1) \}$
its set of bases.  Then both $K[P]$ and
$K[B]$ are Gorenstein.
}
\end{Example}

Let, in general, ${\cal P} \subset \RR^n$
be an integral convex polytope of dimension $n$
and $K[{\cal P}]$ the Ehrhart ring of ${\cal P}$.
A combinatorial criterion
\cite[Corollary (1.2)]{DeNegriHibi}
for $K[{\cal P}]$ to be Gorenstein is stated below.
Let $\delta > 0$ denote the smallest integer
for which
$\delta({\cal P} \setminus \partial{\cal P})
\cap \ZZ^n \neq \emptyset$.  Here
$\delta {\cal P}
= \{ \delta \alpha \: \alpha
\in {\cal P} \}$,
$\partial{\cal P}$ is the boundary of
${\cal P}$ and
${\cal P} \setminus \partial{\cal P}$
is the interior of ${\cal P}$.
Then $K[{\cal P}]$ cannot be Gorenstein unless
$\delta({\cal P} \setminus \partial{\cal P})$
possesses a unique integer vector.
In case that
$\delta({\cal P} \setminus \partial{\cal P})$
possesses a unique integer vector, say
$\alpha_0 \in \ZZ^n$, let
${\cal Q} = \delta{\cal P} - \alpha_0
= \{ \alpha - \alpha_0 \: \alpha
\in \delta{\cal P} \}$.
Thus ${\cal Q} \subset \RR^n$
is an integral convex polytope of dimension
$n$ and the origin of $\RR^n$ is a unique
integer vector belonging to the interior
${\cal Q} \setminus \partial {\cal Q}$
of ${\cal Q}$.
Now, \cite[Corollary (1.2)]{DeNegriHibi} says
that the Ehrhart ring $K[{\cal P}]$
is Gorenstein if and only if the following
condition is satisfied:
If the hyperplane ${\cal H} \subset \RR^n$
determined by a linear equation
$\sum_{i=1}^n a_1 x_i = b$,
where each $a_i$ and $b$ are integers and
where the greatest common divisor of
$a_1, \ldots, a_n, b$ is equal to $1$,
is a supporting hyperplane of ${\cal Q}$
such that ${\cal H} \cap {\cal Q}$
is a facet of ${\cal Q}$, then $b$
is either $1$ or $-1$.

First of all,
we discuss the problem for which discrete
polymatroids $P \subset \ZZ_+^n$ the Ehrhart
ring $K[P]$ is Gorenstein.
In what follows, we will assume that the canonical
basis vectors $\varepsilon_1, \ldots,
\varepsilon_n$ of $\RR^n$ belong to $P$.
In order to apply the above criterion
to the Ehrhart ring $K[P]$, it is required to
know linear equations of supporting
hyperplanes which define facets of
the integral polymatroid
$\conv(P) \subset \RR_+^n$.
Note that, since each $\varepsilon_i \in P$,
the dimension of $\conv(P)$ is equal to $n$.
Let $\rho$ denote the ground set rank function
of $\conv(P)$.
We say that $\emptyset \neq A \subset [n]$
is $\rho$-{\em closed} if any subset
$B \subset [n]$ properly containing $A$
satisfies $\rho(A) < \rho(B)$, and that
$\emptyset \neq A \subset [n]$
is $\rho$-{\em separable} if there exist
two nonempty subsets $A_1$ and $A_2$ of $A$
with $A_1 \cap A_2 = \emptyset$ and
$A_1 \cup A_2 = A$ such that
$\rho(A) = \rho(A_1) + \rho(A_2)$.
For $\emptyset \neq A \subset [n]$
define the hyperplane
${\cal H}_A \subset \RR^n$ by
\[
{\cal H}_A = \{ (x_1, \ldots, x_n) \in \RR^n
\: \sum_{i \in A} x_i = \rho(A) \}.
\]
In addition, for each $i \in [n]$
define the hyperplane
${\cal H}^{(i)} \subset \RR^n$ by
\[
{\cal H}^{(i)} = \{ (x_1, \ldots, x_n) \in \RR^n
\: x_i = 0 \}.
\]
Edmonds \cite{Edmonds} says

\begin{Proposition}
\label{Edmondsfacet}
The facets of $\conv(P) \subset \RR_+^n$ are all
${\cal H}^{(i)} \cap \conv(P)$, where
$i$ ranges over all integers belonging to $[n]$,
and all ${\cal H}_A \cap \conv(P)$, where $A$
ranges over all $\rho$-closed and $\rho$-inseparable
subsets of $[n]$.
\end{Proposition}

We are now in the position to state
a characterization for the Ehrhart ring $K[P]$
of a discrete polymatroid $P \subset \ZZ_+^n$
to be Gorenstein.  For $A \subset [n]$
write $|A|$ for the cardinality of $A$.

\begin{Theorem}
\label{EhrhartGorenstein}
Let $P \subset \ZZ_+^n$ be a discrete polymatroid
and suppose that the canonical basis vectors
$\varepsilon_1, \ldots, \varepsilon_n$ of $\RR^n$
belong to $P$.  Let $\rho$ denote the ground set
rank function of the integral polymatroid
${\cal P} = \conv(P) \subset \RR^n$.
Then the Ehrhart ring
$K[P]$ of $P$ is Gorenstein if and only if
there exists an integer $\delta \geq 1$ such that
\[
\rho(A) = \frac{1}{\delta}(|A| + 1)
\]
for all $\rho$-closed and $\rho$-inseparable
subsets $A$ of $[n]$.
\end{Theorem}

\begin{pf}
Let $\delta \geq 1$ denote the smallest integer
for which
$\delta({\cal P} \setminus \partial {\cal P})
\cap \ZZ^n \neq \emptyset$.  Since
$\varepsilon_1, \ldots, \varepsilon_n$
belong to $P$, it follows that
$\delta({\cal P} \setminus \partial {\cal P})
\cap \ZZ^n \neq \emptyset$
if and only if $(1, \ldots, 1) \in
\delta({\cal P} \setminus \partial {\cal P})$.
By virtue of Proposition \ref{Edmondsfacet}
the equations of the supporting hyperplanes
which define the facets of
$\delta {\cal P} - (1, \ldots, 1)$
are all $x_i = -1$, $1 \leq i \leq n$, and
all $\sum_{i \in A} x_i = \delta \rho(A) - |A|$,
where $A$ ranges over all $\rho$-closed and
$\rho$-inseparable subsets of $[n]$.
Since $(1, \ldots, 1) \in \delta {\cal P}$,
one has $|A| \leq \delta \rho(A)$.
Thus if $K[P]$ is Gorenstein, then
$\delta \rho(A) - |A| = 1$
for all $\rho$-closed and
$\rho$-inseparable subsets of $[n]$.
Conversely,
if there is an integer $\delta \geq 1$ such that
$\delta \rho(A) - |A| = 1$
for all $\rho$-closed and
$\rho$-inseparable subsets of $[n]$,
then $(1, \ldots, 1)$ is a unique integer vector
belonging to
$\delta({\cal P} \setminus \partial {\cal P})$
and $K[P]$ is Gorenstein.
\end{pf}

\begin{Example}
{\em
(a) Let $P_n^{(d)} \subset \ZZ_+^n$ be the
discrete polymatroid consisting of all integer
vectors $u \in \ZZ^n_+$ with $|u| \leq d$
and $B_n^{(d)}$ the set of bases of $P_n^{(d)}$.
Let $\rho$ denote the ground set rank function
of the integral polymatroid
$\conv(P) \subset \RR_+^n$.
Then $\rho(A) = d$ for all
$\emptyset \neq A \subset [n]$.  Thus $[n]$
is the only $\rho$-closed and
$\rho$-inseparable subset of $[n]$.
Hence the Ehrhart ring $K[P_n^{(d)}]$
is Gorenstein if and only if
$d$ divides $n + 1$.  On the other hand,
the base ring $K[B_n^{(d)}]$ is the
Veronese subring $K[t_1, \ldots, t_n]^{(d)}$.
Thus
$K[B_n^{(d)}]$
is Gorenstein if and only if
$d$ divides $n$.  (Note, in fact, that
$K[P_n^{(d)}]$ is just the Veronese subring
$K[x_1, \ldots, x_n, s]^{(d)}$.)

(b) Fix $\delta \geq 1$ which divides $n + 1$.
Let $\rho\: 2^{[n]} \to \RR_+$ denote
the nondecreasing function defined by
\[
\rho(A) = \lceil
( \max(A) + 1 ) / \delta \rceil
\]
for all $\emptyset \neq A \subset [n]$
together with $\rho(\emptyset) = 0$,
where $\max(A)$ is the biggest integer
belonging to $A$ and where $\lceil x \rceil$
is the smallest integer $\geq x$.
Then $\rho$ is submodular.
(In fact, if $A, B \subset [n]$ are nonempty
with $\max(B) \leq \max(A)$, then
$\max(A \cup B) = \max(A)$ and
$\max(A \cap B) \leq \max(B)$.
Thus $\rho(A \cup B) = \rho(A)$ and
$\rho(A \cap B) \leq \rho(B)$.)
If $\emptyset \neq A \subset [n]$
is $\rho$-closed, then
$A = [r]$ for some $1 \leq r \leq n$
such that $\delta$ divides $r + 1$.
Moreover, if $\delta$ divides $r + 1$,
then $A = [r]$ satisfies
$\rho(A) = (|A| + 1) / \delta$.
Let $d = (n + 1) / \delta$.  Let
$P \subset \ZZ_+^n$ be the discrete
polymatroid of rank $d$
arising from $\rho$, i.e.,
$P$ consists of
all integer vectors $u \in \ZZ_+^n$
such that $u(A) \leq \rho(A)$
for all $\emptyset \neq A \subset [n]$.
More precisely,
$u = (u_1, \ldots, u_n) \in \ZZ_+^n$
belongs to $P$ if and only if
\[
u_1 + u_2 + \cdots + u_{i \delta - 1}
\leq i, \, \, \, \, \,
i = 1, 2, \cdots, d.
\]
It follows from Theorem \ref{EhrhartGorenstein}
that the Ehrhart ring $K[P]$ is Gorenstein.

We also discuss the base ring
of this discrete polymatroid.
Let $\delta \geq 2$.  Then
$u \in \ZZ_+^n$ is a base of $P$
if and only if
$u$ is of the form
\[
u = \varepsilon_{i_1} + \varepsilon_{i_2}
+ \cdots + \varepsilon_{i_d},
\]
where $1 \leq i_1 \leq \delta - 1$ and
where
$\delta (j - 1) \leq i_j \leq \delta j - 1$
for $j = 2, \ldots, d$.
Thus the base ring $K[B]$ of $P$
is the Segre product
\[
K[t_1, \ldots, t_{\delta-1}] \sharp
K[t_\delta, \ldots, t_{2\delta-1}] \sharp
\cdots \sharp
K[t_{(d-1)\delta}, \ldots, t_{d \delta-1}]
\]
of the polynomial ring in $\delta - 1$ variables
and $d - 1$ copies of
the polynomial ring in $\delta$ variables.
Thus $K[B]$ is not Gorenstein
unless $\delta = n + 1$.
If $\delta = 1$, then $(2, 1, \ldots, 1)$
is a unique base of $P$ and the base ring
$K[B] = K[t_1^2 t_2 \cdots t_n]$
is Gorenstein.

(c) Let $\rho\: 2^{[n]} \to \RR_+$ denote
the nondecreasing function defined by
$\rho(A) = |A| + 1$
for all $\emptyset \neq A \subset [n]$
together with $\rho(\emptyset) = 0$.
Then $\rho$ is submodular and
all nonempty subsets of $[n]$
are $\rho$-closed and $\rho$-inseparable.
Let $P \subset \ZZ_+^n$ be the discrete
polymatroid of rank $n + 1$
arising from $\rho$.
Then the Ehrhart ring $K[P]$ is
Gorenstein.  Moreover, since the set of bases
of $P$ is
$B = \{(2, 1, \ldots, 1),
\ldots, (1, \ldots, 1, 2)\}$,
the base ring $K[B]$ isomorphic to the polynomial
ring in $n$ variables; thus $K[B]$
is Gorenstein.
}
\end{Example}

We now turn to the problem when the base ring of
a discrete polymatroid is Gorenstein.  To obtain
a perfect answer to this problem seems, however,
quite difficult. For example, the discussions
appearing in \cite{DeNegriHibi} for classifying
Gorenstein rings belonging to the class of
algebras of Veronese type,
a distinguished class of discrete polymatroids
(Example \ref{examples1} (c)), is enough complicated.
In what follows we introduce the concept of
``generic'' discrete polymatroids and find
a characterization for the base ring of
a generic discrete polymatroid to be Gorenstein.

Let $P \subset \ZZ_+^n$ be a discrete polymatroid
of rank $d$ and $B = B(P)$ the set of bases of $P$.
We will, as before, assume that the canonical
basis vectors $\varepsilon_1, \ldots, \varepsilon_n$
of $\RR^n$ belong to $P$.
Let ${\cal F} = \conv(B)$, the set of bases of
the integral polymatroid
${\cal P} = \conv(P) \subset \RR_+^n$.
Recall that ${\cal F}$ is a face of ${\cal P}$
with the supporting hyperplane ${\cal H}_{[n]}
\subset \RR^n$, i.e.,
${\cal F} = {\cal H}_{[n]} \cap {\cal P}$.
Let $\rho \: 2^{[n]} \to \RR_+$ denote
the ground set rank function of ${\cal P}$.
Then
\[
{\cal F} = \{ u \in {\cal H}_{[n]} \cap \ZZ_+^n
\: u(A) \leq \rho(A),
\emptyset \neq A \subset [n], A \neq [n]  \}.
\]
Let $\varphi \: {\cal H}_{[n]} \to
\RR^{n-1}$ denote the affine transformation
defined by
\[
\varphi(u_1, \ldots, u_n) = (u_1, \ldots, u_{n-1}).
\]
Thus $\varphi$ is injective and
$\varphi({\cal H}_{[n]} \cap \ZZ^n) = \ZZ^{n-1}$.
Since for all $A \subset [n]$ with $n\in A$ and  $A \neq [n]$ the hyperplane
$\varphi({\cal H}_{A} \cap {\cal H}_{[n]})
\subset \RR^{n-1}$
is determined by the linear equation
$\sum_{i \in [n] \setminus A} x_i
= d - \rho(A)$, it follows that
\[
\varphi({\cal F})
= \{ u \in \RR_+^{n-1} \:
d - \rho([n] \setminus A) \leq u(A) \leq
\rho(A), \emptyset \neq A \subset [n - 1] \}.
\]
We say that $P$ is {\em generic} if
\begin{enumerate}
\item[(G1)]
each base $u$ of $P$ satisfies $u(i) > 0$
for all $1 \leq i \leq n$;
\item[(G2)]
${\cal F} = \conv(B)$ is a facet of
${\cal P} = \conv(P)$;
\item[(G3)]
${\cal F} \cap {\cal H}_{A}$ is a facet of
${\cal F}$ for all $\emptyset \neq A \subset [n]$
with $A \neq [n]$.
\end{enumerate}
It follows that $P$ is generic
if and only if
(i) $\rho$ is strictly increasing,
(ii) $\dim \varphi({\cal F}) = n - 1$,
and (iii)
the facets of $\varphi({\cal F})$ are all
$\{ u \in \varphi({\cal F}) \:
u(A) = \rho(A) \}$
and all
$\{ u \in \varphi({\cal F}) \:
u(A) = d - \rho([n] \setminus A) \}$,
where $A$ ranges over all nonempty subsets
of $[n - 1]$.

\begin{Example}
\label{examplegeneric}
{\em
(a) Let $n = 2$ and let $a_1, a_2 > 0$ be integers.
Let $P \subset \ZZ_+^2$ denote the discrete
polymatroid of rank $d$ consisting of
those $u = (u_1, u_2) \in \ZZ_+^2$ such that
$u_1 \leq a_1$, $u_2 \leq a_2$ and $u_1 + u_2 \leq d$.
Then $P$ is generic if and only if
$a_1 < d$, $a_2 < d$, and $d < a_1 + a_2$.
If $P$ is generic, then
the bases of $P$ are
$(a_1, d - a_1), (a_1 - 1, d - a_1 + 1), \ldots,
(d - a_2, a_2)$.  Thus the base ring $K[B]$
of $P$ is Gorenstein if and only if
either $a_1 + a_2 = d + 1$ or $a_1 + a_2 = d + 2$.

(b) Let $n = 3$.  Let $P \subset \ZZ_+^3$
be a discrete polymatroid of rank $d$ with $B$
its set of bases,
and $\rho$ the ground set rank function
of the integral polymatroid
$\conv(P) \subset \RR_+^3$.
Then $\varphi({\cal F}) \subset \ZZ_+^2$,
where ${\cal F} = \conv(B)$, consists of
those $u = (u_1, u_2) \in \RR_+^2$ such that
\begin{eqnarray*}
d - \rho(\{ 2, 3 \})
\leq & u_1
\, \, \, \, \, \, \, \, \, \, \, \, \,
& \leq \rho(\{ 1 \}),
\\
d - \rho(\{ 1, 3 \})
\leq & \, \, \, \, \, \, \, \, \, \, \, \, \,
u_2
& \leq \rho(\{ 2 \}),
\\
d - \rho(\{ 3 \})
\leq & u_1 + u_2
& \leq \rho(\{ 1, 2 \}).
\end{eqnarray*}
Hence $P$ is generic if and only if
\begin{eqnarray*}
0 < \rho(\{ i \}) < \rho(\{ i, j \}) < d,
& &
1 \leq i \neq j \leq 3, \\
\rho(\{ i \}) + \rho(\{ j \}) > \rho(\{ i, j \}),
& &
1 \leq i < j \leq 3, \\
\rho(\{ i, j \}) + \rho(\{ j, k  \}) >
d + \rho(\{ j \}),
& &
\{ i, j, k \} = [3].
\end{eqnarray*}
Moreover, if $P$ is generic, then the base ring
$K[B]$ is Gorenstein if and only if
\[
\rho(\{ i \}) + \rho(\{ j, k \}) = d + 2,
\, \, \, \, \, \{ i, j, k \} = [3].
\]
}
\end{Example}

\begin{Theorem}
\label{genericGorenstein}
{\em (a)}
Let $n \geq 3$.
Let $P \subset \ZZ_+^n$
be a discrete polymatroid of rank $d$
and suppose that the canonical basis vectors
$\varepsilon_1, \ldots, \varepsilon_n$
of $\RR^n$ belong to $P$.
Let $\rho \: 2^{[n]} \to \RR_+$
denote the ground set rank function of
the integral polymatroid
$\conv(P) \subset \RR_+^n$.
If $P$ is generic and if the base ring $K[B]$
of $P$ is Gorenstein, then
there is a vector
$\alpha = (\alpha_1, \ldots, \alpha_{n-1})
\in \ZZ_+^{n-1}$ with each $\alpha_i > 1$
and with $d > |\alpha| + 1$ such that
\[
\rho(A)=\left\{\begin{array}{ll}
\alpha(A) + 1, &\text{if $\quad\emptyset \neq A \subset [n - 1]$}, \\
d - \alpha([n] \setminus A) + 1,& \text{if $\quad n \in A \neq [n].$}
\end{array}
\right.
\]
{\em (b)}
Conversely, given
$\alpha = (\alpha_1, \ldots, \alpha_{n-1})
\in \ZZ_+^{n-1}$, where $n \geq 3$,
with each $\alpha_i > 1$
and $d \in \ZZ$ with $d > |\alpha| + 1$,
define the function
$\rho \: 2^{[n]} \to \RR_+$
by (a) together with
$\rho(\emptyset) = 0$ and $\rho([n]) = d$.
Then
\begin{enumerate}
\item[(i)]
$\rho$ is strictly increasing and submodular;
\item[(ii)]
the discrete polymatroid
$P = \{ u \in \ZZ_+^n \:
u(A) \leq \rho(A),
\emptyset \neq A \subset [n] \}
\subset \ZZ_+^n$
arising from $\rho$ is generic;
\item[(iii)]
the base ring $K[B]$ of $P$
is Gorenstein.
\end{enumerate}
\end{Theorem}

\begin{pf}
(a)
Suppose that a discrete polymatroid
$P \subset \ZZ_+^n$ of rank $d$
is generic and that the base ring $K[B]$
of $P$ is Gorenstein.
Let ${\cal F} = \conv(B)$.
Since $K[B]$ is Gorenstein, there is
an integer $\delta \geq 1$ such that
\[
\delta ( \rho(A) - ( d - \rho([n] \setminus A) ) )
= 2
\]
for all $\emptyset \neq A \subset [n - 1]$.  Hence
either $\delta = 1$ or $\delta = 2$.

If $\delta = 2$, then
$(2 \rho(\{ 1 \}) - 1, \ldots,
2 \rho(\{ n - 1 \}) - 1) \in \ZZ_+^{n - 1}$
must be a unique integer vector belonging to
the interior of $2 \varphi({\cal F})$.
Since $K[B]$ is Gorenstein it follows that 
\[
\sum_{i \in A} ( 2 \rho(\{ i \}) - 1 )
= 2 \rho(A) - 1,
\, \, \, \, \,
\emptyset \neq A \subset [n - 1].
\]
Thus
\[
\rho(A) =
\sum_{i \in A} \rho(\{ i \})
- \frac{1}{2} (|A| - 1),
\, \, \, \, \,
\emptyset \neq A \subset [n - 1].
\]
Since $n \geq 3$, it follows that
$\rho(\{1, 2\}) \not\in \ZZ$, a contradiction.

Now let $\delta = 1$ and set
\[
\alpha_i = \rho(\{ i \}) - 1 =
d - \rho([n] \setminus \{ i \}) + 1 > 1,
\, \, \, \, \,
1 \leq i \leq n - 1.
\]
Then $\alpha = (\alpha_1, \ldots, \alpha_{n-1})
\in \ZZ_+^{n-1}$
is a unique integer vector belonging
to the interior of $\varphi({\cal F})
\subset \RR_+^{n-1}$ and
$\varphi({\cal F}) - \alpha$ consists of
those $u = (u_1, \ldots, u_{n-1}) \in
\RR^{n-1}$ such that
\[
d - \rho([n] \setminus A) - \alpha(A)
\leq \sum_{i \in A} u_i \leq
\rho(A) - \alpha(A),
\, \, \, \, \,
\emptyset \neq A \subset [n - 1].
\]
Since $P$ is generic, the desired equality
on $\rho$ follows immediately.
Moreover,
since
$\rho([n-1]) = |\alpha| + 1 < \rho([n]) = d$,
one has $d > |\alpha| + 1$, as required.

(b) Since each $\alpha_i > 1$ and since
$d > |\alpha| + 1$,
it follows that
$0 < \rho(A) < \rho([n]) = d$
for all $\emptyset \neq A \subset [n]$
with $A \neq [n]$.
Moreover, $\rho(\{ i \}) > 2$
for all $1 \leq i \leq n$.
If $\emptyset \neq A \subset B \subset [n]$
with $n \not\in A$ and $n \in B$, then
$\rho(B) - \rho(A) = d -
(\alpha([n] \setminus B) + \alpha(A)) >
d - |\alpha| > 1$.
Hence $\rho$ is strictly increasing.

To see why $\rho$ is submodular,
we distinguish three cases as follows.
First, if $A, B \subset [n - 1]$
with $A \neq \emptyset$ and
$B \neq \emptyset$, then
\begin{eqnarray*}
\rho(A) + \rho(B) & = & \alpha(A) +
\alpha(B) + 2 \\
& = & \alpha(A \cup B) + \alpha(A \cap B)
+ 2 \\ & = &
\rho(A \cup B) + \rho(A \cap B)
\end{eqnarray*}
unless $A \cap B \neq \emptyset$.
Second, if $A, B \subset [n]$ with
$n \in A \neq [n]$ and
$n \in B \neq [n]$, then
\begin{eqnarray*}
\rho(A) + \rho(B) & = &
2d - (\alpha([n] \setminus A)
+ \alpha([n] \setminus B)) + 2 \\
& = &
2d - (\alpha([n] \setminus (A \cap B))
+ \alpha([n] \setminus (A \cup B)))
+ 2 \\ & = & \rho(A\cup B) + \rho(A \cap B)
\end{eqnarray*}
unless $A \cup B \neq [n]$.
Third, if $n \in A$ and $B \subset [n - 1]$,
then assuming $A \cap B \neq \emptyset$
and $A \cup B \neq [n]$ one has
\begin{eqnarray*}
\rho(A) + \rho(B) & = &
d - \alpha([n] \setminus A) + 1
+ \alpha(B) + 1 \\
& = &
d - \alpha(([n] \setminus A) \setminus B)
+ 1
+ \alpha(B \setminus ([n] \setminus A))
+ 1 \\
& = &
d - \alpha([n] \setminus (A \cup B)) + 1
+ \alpha(A \cap B) + 1 \\
& = & \rho(A\cup B) + \rho(A \cap B),
\end{eqnarray*}
as desired.

Let ${\cal F} = \conv(B)$.  Since
\[
\varphi({\cal F}) =
\{ u \in \RR_+^{n-1} \:
\alpha(A) - 1 \leq u(A) \leq \alpha(A) + 1,
\emptyset \neq A \subset [n - 1] \},
\]
it follows that
$\varphi({\cal F}) - \alpha
\subset \RR^{n-1}$
consists of those
$u = (u_1, \ldots, u_{n-1}) \in \RR^{n-1}$
such that
\[
-1 \leq u_{i_1} + \cdots +  u_{i_k} \leq 1,
\, \, \, \, \,
1 \leq i_1 < \cdots < i_k \leq n - 1.
\]
Hence
$(\varphi({\cal F}) - \alpha) \cap \ZZ^{n-1}$
consists of those
$v = (v_1, \ldots, v_{n-1}) \in \ZZ^{n-1}$
such that
$-1 \leq v_i \leq 1$ for all $1 \leq i \leq n - 1$,
$| \{ i \: v_i = 1 \}| \leq 1$, and
$| \{ i \: v_i = - 1 \}| \leq 1$.
In particular the canonical basis vectors
$\varepsilon_1, \ldots, \varepsilon_{n-1}$
of $\RR^{n-1}$ belong to
$\varphi({\cal F}) - \alpha$.
Thus ${\cal F}$ is a facet of
${\cal P} = \conv(P)$.
For $1 \leq i_1 < \cdots < i_k \leq n - 1$
write ${\cal H}_{i_1 \cdots i_k} \subset \RR^{n-1}$
(resp.
${\cal H}'_{i_1 \cdots i_k} \subset \RR^{n-1}$)
for the supporting hyperplane of ${\cal F}$
determined by the linear equation
$x_{i_1} + \cdots + x_{i_k} = 1$
(resp. $x_{i_1} + \cdots + x_{i_k} = - 1$).
Then the vectors
$\varepsilon_{i_1}, \ldots, \varepsilon_{i_k}$
(resp.
$- \varepsilon_{i_1}, \ldots, - \varepsilon_{i_k}$)
and
$\varepsilon_{i_1} - \varepsilon_j$
(resp.
$- \varepsilon_{i_1} + \varepsilon_j$),
$j \in [n-1] \setminus \{ i_1, \ldots, i_k \}$,
belong to the face
$(\varphi({\cal F}) - \alpha) \cap
{\cal H}_{i_1 \cdots i_k}$
(resp.
$(\varphi({\cal F}) - \alpha) \cap
{\cal H}'_{i_1 \cdots i_k}$)
of $\varphi({\cal F}) - \alpha$.
Thus
$(\varphi({\cal F}) - \alpha) \cap
{\cal H}_{i_1 \cdots i_k}$
(resp.
$(\varphi({\cal F}) - \alpha) \cap
{\cal H}'_{i_1 \cdots i_k}$)
is a facet of $\varphi({\cal F}) - \alpha$.
Hence $P$ is generic.
Moreover, since the Ehrhart ring
$K[\varphi({\cal F}) - \alpha]$
is Gorenstein,
the base ring $K[B]$
$(\cong K[\varphi({\cal F}) - \alpha])$
is Gorenstein, as desired.
\end{pf}

\section{Constructions of discrete polymatroids}
In this section two simple
techniques to construct discrete polymatroids
will be studied.  The first one shows that
a nondecreasing and submodular function
defined in a sublattice of $2^{[n]}$
produces a discrete polymatroid.
The second one yields the concept of
transversal polymatroids.

A {\em sublattice} of $2^{[n]}$ is
a collection ${\cal L}$ of subsets
of $[n]$ with $\emptyset \in {\cal L}$
and $[n] \in {\cal L}$ such that
for all $A$ and $B$ belonging to ${\cal L}$
both $A \cap B$ and $A \cup B$ belong
to ${\cal L}$.
A function $\mu \: {\cal L} \to \RR_+$
is called {\em submodular} if
\[
\mu(A) + \mu(B) \geq
\mu(A \cup B) + \mu(A \cap B),
\, \, \, \, \,
A, B \in {\cal L}.
\]
Even though Theorem \ref{sublattice} below
is a direct consequence of
\cite[Theorem 1, p. 342]{Welsh}
together with Theorem \ref{firstmaintheorem},
for the sake of completeness
we will give an easy proof based on
\cite[Theorem 2, p. 345]{Welsh}.

\begin{Theorem}
\label{sublattice}
Let ${\cal L}$ be a sublattice of $2^{[n]}$
and $\mu \: {\cal L} \to \RR_+$
an integer valued nondecreasing and
submodular function
with $\mu(\emptyset) = 0$.
Then
\[
P_{({\cal L},\mu)}
=
\{ u \in \ZZ_+^n \:
u(A) \leq \mu(A),\quad A \in {\cal L} \}.
\]
is a discrete polymatroid.
\end{Theorem}

\begin{pf}
Let $\rho \: 2^{[n]} \to \RR_+$
be the nondecreasing function define by
\[
\rho(X) = \min \{ \mu(A) \:
A \supset X,\quad A \in {\cal L} \}
\]
together with $\rho(\emptyset) = 0$.
Then $\rho$ is submodular.
Let $P_{\rho} \subset \ZZ_+^{n}$
denote the discrete polymatroid
arising from $\rho$.
Since $\emptyset \neq X \subset [n]$
with $X \not\in {\cal L}$ cannot be
$\rho$-closed,
it follows from Proposition
\ref{Edmondsfacet}
that
\[
P_{\rho}
= \{ u \in \ZZ_+^{n} \:
u(A) \leq \rho(A),\quad
A \in {\cal L} \}.
\]
Since $\rho(A) = \mu(A)$
for all $A \in {\cal L}$,
one has
$P_{\rho} = P_{({\cal L},\mu)}$,
thus $P_{({\cal L},\mu)}$
is a discrete polymatroid,
as desired.
\end{pf}

\begin{Example}
\label{principalBorel}
{\em
Let ${\cal L}$ be a chain of length $n$
of $2^{[n]}$, say
$${\cal L} = \{ \emptyset,
\{ n \},
\{ n - 1 , n \}, \ldots,
\{ 1, \ldots, n \}
\} \subset 2^{[n]}.$$
Given nonnegative integers
$a_1, \ldots, a_n$, define
$\mu \: {\cal L} \to \RR_+$
by
$$\mu(\{ i, i + 1, \ldots, n \})
= a_i + a_{i+1} + \cdots + a_n,
\, \, \, \, \,
1 \leq i \leq n$$
together with $\mu(\emptyset) = 0$.
Then the discrete polymatroid
$P_{({\cal L},\mu)} \subset \ZZ_+^n$
is
\[
P_{({\cal L},\mu)} =
\{ u = (u_1, \ldots, u_n) \in \ZZ_+^n \:
\sum_{j=i}^{n} u_j \leq
\sum_{j=i}^{n} a_j,
\, \, \, 1 \leq i \leq n \}.
\]
This example will be discussed again in
Example \ref{transexample}.
}
\end{Example}

Let ${\cal A} =
(A_1, \ldots, A_d)$
be a family of nonempty subsets of $[n]$.
It is {\em not} required that
$A_i \neq A_j$ if $i \neq j$.
Let
\[
B_{\cal A} =
\{
\varepsilon_{i_1} + \cdots +
\varepsilon_{i_d}
\:
i_k \in A_k,\quad  1 \leq k \leq d
\}
\subset \ZZ_+^n
\]
and define the integer valued nondecreasing
function
$\rho_{\cal A} \: 2^{[n]} \to \RR_+$
by setting
\[
\rho_{\cal A}(X) =
| \{ k \: A_k \cap X \neq \emptyset \} |,
\, \, \, \, \,
X \subset [n].
\]

\begin{Theorem}
\label{transversal}
The function $\rho_{\cal A}$ is submodular and
$B_{\cal A}$ is the set of bases of the
discrete polymatroid $P_{\cal A} \subset \ZZ_+^n$
arising from $\rho_{\cal A}$.
\end{Theorem}

\begin{pf}
For each $1 \leq k \leq d$ define the function
$\rho_{k} \: 2^{[n]} \to \RR_+$
by $\rho_k(X) = 1$
if $A_k \cap X \neq \emptyset$ and
$\rho_k(X) = 0$ if $A_k \cap X = \emptyset$.
Then $\rho_k$ is nondecreasing and submodular.
Write $P_k \subset \ZZ_+^n$
for the discrete polymatroid arising from
$\rho_k$.  Then Proposition \ref{sum} guarantees
that
\[
P = \{ x \in \ZZ_+^n \:
x = \sum_{k=1}^d x_k,
\, \, \, x_k \in P_k \}
\]
is a discrete polymatroid of rank $d$
and the ground set
rank function of the integral polymatroid
$\conv(P) \subset \RR_+^n$
is $\rho = \sum_{k=1}^d \rho_k$.  It is clear that
\[
\rho(X) =
| \{ k \: A_k \cap X \neq \emptyset \} |,
\, \, \, \, \,
X \subset [n]
\]
and the set of bases of $P$ coincides with
$B_{\cal A}$, as desired.
\end{pf}

We say that the discrete polymatroid
$P_{\cal A} \subset \ZZ_+^n$
is the {\em transversal polymatroid}
presented by ${\cal A}$.

\begin{Example}
\label{transexample}
{\em
Let $r_1, \ldots, r_d \in [n]$ and
set $A_k = [r_k] = \{ 1, \ldots, r_k \}$,
$1 \leq k \leq d$.
Let $\min(X)$ denote the smallest integer
belonging to $X$, where
$\emptyset \neq X \subset [n]$.  If ${\cal A} =
(A_1, \ldots, A_d)$, then
\[
\rho_{\cal A}(X) = \rho_{\cal A}(\{\min(X)\})
= | \{ k \: \min(X) \leq r_k \}|.
\]
If $\emptyset \neq X \subset [n]$
is $\rho_{\cal A}$-closed, then
$X = \{ \min(X), \min(X) + 1, \ldots, n \}$.
Let
$$a_i = |\{ k \: r_k = i \}|,
\, \, \, \, \,
1 \leq i \leq n.$$
Thus
$$\rho_{\cal A}(\{i, i+1, \ldots, n\})
= a_i + a_{i+1} + \cdots + a_n,
\, \, \, \, \,
1 \leq i \leq n.$$
The transversal polymatroid
$P_{\cal A} \subset \ZZ_+^n$
presented by ${\cal A}$ is
\[
P_{\cal A} =
\{ u = (u_1, \ldots, u_n) \in \ZZ_+^n \:
\sum_{j=i}^{n} u_j \leq
\sum_{j=i}^{n} a_j,
\, \, \, 1 \leq i \leq n \}.
\]
Thus $P_{\cal A}$ coincides with the
discrete polymatroid $P_{({\cal L}, \mu)}$
in Example \ref{principalBorel}.
If
$P_i \subset \ZZ_+^n$
is the discrete polymatroid
with $B_i = \{ \varepsilon_1, \ldots,
\varepsilon_i \}$
its set of bases, then
$P_{\cal A} = a_1 P_1 \vee \cdots \vee
a_n P_n$.  Moreover, $B_{\cal A}$
is equal to the strongly stable principal set
with $(a_1, \ldots, a_n)$ its Borel generator.
}
\end{Example}

\begin{Example}
\label{nontransversal}
{\em
Let $P \subset \ZZ_+^4$ denote the discrete
polymatroid of rank $3$ consisting of those
$u = (u_1, u_2, u_3, u_4) \in \ZZ_+^4$
with $u_i \leq 2$ for $1 \leq i \leq 4$
and with $|u| \leq 3$.  Then
$P$ is {\em not} transversal.
Suppose, on the contrary, that
$P$ is a transversal polymatroid
presented by ${\cal A} = (A_1, A_2, A_3)$
with each $A_k \subset [4]$.
Since
$(2,1,0,0), (2,0,1,0), (2,0,0,1) \in P$
and $(3,0,0,0) \not\in P$,
we assume that $1 \in A_1$, $1 \in A_2$
and $A_3 = \{ 2, 3, 4 \}$.
Since
$(1,2,0,0), (0,2,1,0), (0,2,0,1) \in P$
and $(0,3,0,0) \not\in P$,
we assume that $2 \in A_1$
and $A_2 = \{ 1, 3, 4 \}$.
Since $(0,0,2,1) \in P$
and $(0,0,3,0) \not\in P$,
one has $4 \in A_1$.  Hence
$(0,0,0,3) \in P$, a contradiction.
}
\end{Example}

The discussions in Example \ref{transexample}
enables us to classify all Gorenstein strongly
stable principal sets.

\begin{Proposition}
\label{principalBorelGorenstein}
Let $a = (a_1, \ldots, a_n) \in \ZZ_+^n$
with $a_n \geq 1$ and write $B_a$
for the strongly stable principal set
with $(a_1, \ldots, a_n)$ its Borel generator.
Then the base ring $K[B_a]$ is Gorenstein
if and only if
$a_2 + \cdots + a_n$ divides $n$ and
\[
\frac{n - i + 2}{a_i + a_{i+1} + \cdots + a_n}
=
\frac{n}{a_2 + \cdots + a_n}
\]
for all $3 \leq i \leq n$ with
$a_{i-1} \neq 0$.
\end{Proposition}

\begin{pf}
The base ring
$K[B_a]$ is isomorphic to the Ehrhart ring
of the integral convex polytope
${\cal P} \in \RR^{n-1}$ consisting of those
$u = (u_2, u_3, \ldots, u_n)$ such that
$u_i \geq 0$, $2 \leq i \leq n$, and
$\sum_{j=i}^{n} u_j \leq
\sum_{j=i}^{n} a_j$ for $i = 2$
and for all $3 \leq i \leq n$
with $a_{i-1} \neq 0$.
Since all of these inequalities define
the facets of ${\cal P}$, the desired result
follows immediately.
Note that
${\cal P}$
is an integral polymatroid
whose ground set rank function
$\rho$ is given by
$\rho(X) = \sum_{j = \min(X)}^n
a_j$ for
$\emptyset \neq X \subset \{ 2, \ldots, n \}$.
\end{pf}

\begin{Example}
\label{easyexample}
{\em
(a) If $a = (0, 1, \ldots, 1 , 2)
\in \ZZ_+^n$, then $K[B_a]$ is Gorenstein.

(b) If $a = (0, 1, 0, 2, 0, 3)$, then
$K[B_a]$ is Gorenstein.

(c) If $a = (0, \ldots, 0 , a_n)
\in \ZZ_+^n$ with $a_n \geq 1$, then
$K[B_a]$ is Gorenstein
if and only if $a_n$ divides $n$.
}
\end{Example}

It would, of course, be of interest
to classify all transversal polymatroids
with Gorenstein base rings.


\begin{thebibliography}{99}

\bibitem{BruHer} W.\ Bruns and J.\ Herzog, ``Cohen--Macaulay
rings,''  Revised Edition, Cambridge University Press, 
Cambridge, 1996. 

\bibitem{ConHer} A.\ Conca and J.\ Herzog,
Castelnuovo--Mumford regularity of products of ideals,
preprint (2001).



\bibitem{DeNegri} E.\ De Negri, 
Toric rings generated by special stable sets of monomials, 
{\em Math. Nachr.} {\bf 203}
(1999), 31 -- 45.

\bibitem{DeNegriHibi} E.\ De Negri and T.\ Hibi,
Gorenstein algebras of Veronese type,
{\em J. Algebra} {\bf 193} (1997), 629 -- 639.

\bibitem{Edmonds} J.\ Edmonds, Submodular functions,
matroids, and certain polyhedra, {\em in}
``Combinatorial Structures and Their Applications''
(R.\ Guy, H.\ Hanani, N.\ Sauer and J.\ Schonheim, Eds.)
Gordon and Breach, New York, 1970, pp. 69 -- 87.

\bibitem{HerTak} J.\ Herzog and Y.\ Takayama,
Resolutions by mapping cones, to appear.

\bibitem{Hibi} T.\ Hibi, ``Algebraic Combinatorics on
Convex Polytopes,'' Carslaw, Glebe, N.S.W., Australia, 1992.

\bibitem{OhHerHi} H.\ Ohsugi, J.\ Herzog and T.\ Hibi,
Combinatorial pure subrings,
{\em Osaka J.\ Math.} {\bf 37} (2000), 745 -- 757.

\bibitem{Oxley} J. G.\ Oxley, ``Matroid Theory,''
Oxford University Press, Oxford, New York, 1992.

\bibitem{Sturmfels} B.\ Sturmfels, 
``Gr\"obner Bases and Convex Polytopes,'' 
Amer. Math. Soc., Providence, RI, 1995.

\bibitem{Welsh} D. J. A.\ Welsh, ``Matroid Theory,''
Academic Press, London, New York, 1976.  

\bibitem{White} N.\ White, A unique exchange property for bases,
{\em Linear Algebra Appl.} {\bf 31} (1980), 81 -- 91. 




\end{thebibliography}
\end{document}